\documentclass[11pt,pdftex]{article}
\usepackage{latexsym,amssymb}
\usepackage{psfrag}
\usepackage{amsmath}
\usepackage{color}
\usepackage{subfig}
\usepackage{tikz}
\usetikzlibrary{intersections, calc, angles, arrows.meta}
\usetikzlibrary{patterns}

\newtheorem{Theorem}{\bf Theorem}[section]
\newtheorem{Lemma}{\bf Lemma}[section]
\newtheorem{Proposition}{\bf Proposition}[section]
\newtheorem{Corollary}{\bf Corollary}[section]
\newtheorem{Remark}{\bf Remark}[section]
\newtheorem{Example}{\bf Example}[section]
\newtheorem{Definition}{\bf Definition}[section]

\newenvironment{theorem}{\begin{Theorem}$\!\!\!$}{\end{Theorem}}
\newenvironment{lemma}{\begin{Lemma}$\!\!\!$}{\end{Lemma}}
\newenvironment{proposition}{\begin{Proposition}$\!\!\!$}{\end{Proposition}}

\newenvironment{remark}{\begin{Remark}$\!\!\!$}{\end{Remark}}

\newenvironment{definition}{\begin{Definition}$\!\!\!$}{\end{Definition}}

\def\XXint#1#2#3{{\setbox0=\hbox{$#1{#2#3}{\int}$}
\vcenter{\hbox{$#2#3$}}\kern-.5\wd0}}

\numberwithin{equation}{section}

\allowdisplaybreaks[4]

\begin{document}

\title{
Optimal singularities of initial functions\\
for solvability of a semilinear parabolic system}
\author{
Yohei Fujishima and Kazuhiro Ishige
}
\date{}
\maketitle
\begin{abstract}
Let $(u,v)$ be a nonnegative solution to the semilinear parabolic system
$$
\mbox{(P)}
\qquad
\left\{
\begin{array}{ll}
\partial_t u=D_1\Delta u+v^p,\quad & x\in{\bf R}^N,\,\,\,t>0,\vspace{3pt}\\
\partial_t v=D_2\Delta v+u^q,\quad &  x\in{\bf R}^N,\,\,\,t>0,\vspace{3pt}\\
(u(\cdot,0),v(\cdot,0))=(\mu,\nu),\quad & x\in{\bf R}^N,
\end{array}
\right.
$$
where $D_1$, $D_2>0$, $0<p\le q$ with $pq>1$ and $(\mu,\nu)$ is a pair of
nonnegative Radon measures or
nonnegative measurable functions in ${\bf R}^N$.
In this paper we study sufficient conditions on the initial data for the solvability of problem~(P)
and clarify optimal singularities of the initial functions for the solvability.
\end{abstract}
\vspace{50pt}
\noindent Addresses:

\smallskip
\noindent Y.~F.:  Department of Mathematical and Systems Engineering, Faculty of Engineering, Shizuoka University,
3-5-1, Johoku, Hamamatsu, 432-8561, Japan.\\
\noindent
E-mail: {\tt fujishima@shizuoka.ac.jp}\\

\smallskip
\noindent K.~I.:  Graduate School of Mathematical Sciences, The University of Tokyo,
3-8-1 Komaba, Meguro-ku, Tokyo 153-8914, Japan.\\
\noindent
E-mail: {\tt ishige@ms.u-tokyo.ac.jp}
\newpage
\section{Introduction}
We are concerned with the Cauchy problem for a semilinear parabolic system
\begin{equation}
\tag{P}
\label{eq:P}
\left\{
\begin{array}{ll}
\partial_t u=D_1\Delta u+v^p & \quad\mbox{in}\quad{\bf R}^N\times(0,T),\vspace{3pt}\\
\partial_t v=D_2\Delta v+u^q & \quad\mbox{in}\quad{\bf R}^N\times(0,T),\vspace{3pt}\\
u,v\ge 0 & \quad\mbox{in}\quad{\bf R}^N\times(0,T),\vspace{3pt}\\
(u(\cdot,0),v(\cdot,0))=(\mu,\nu) & \quad\mbox{in}\quad{\bf R}^N,
\end{array}
\right.
\end{equation}
where $N\ge 1$, $0<T\le\infty$, $D_1$, $D_2>0$, $0<p\le q$ with $pq>1$
and $(\mu,\nu)$ is a pair of (nonnegative) Radon measures or measurable functions in ${\bf R}^N$.
Problem~(P) is one of the simplest parabolic systems
and it is an example of reaction-diffusion systems describing heat propagation
in a two component combustible mixture.
It has been studied extensively in many papers from various points of view,
see e.g. \cite{AHV, EH01, EH02, FIBu, FIM, IKS, QS} and \cite[Section~32]{QSBook}.
The following properties have been already proved for the case of $D_1=D_2$.
\begin{itemize}
  \item[{\rm (1)}]
  Let $p\ge1$ and $r_1,r_2\in(1,\infty)$.
  Assume
  $$
  \max\{P(r_1,r_2),Q(r_1,r_2)\}\le 2,
  $$
  where
  $$
  P(r_1,r_2):=N\left(\frac{p}{r_2}-\frac{1}{r_1}\right),
  \qquad
  Q(r_1,r_2):=N\left(\frac{q}{r_1}-\frac{1}{r_2}\right).
 $$
  Then, for any $(\mu,\nu)\in L^{r_1,\infty}({\bf R}^N)\times L^{r_2,\infty}({\bf R}^N)$,
  problem~(P) possesses a solution for some $T>0$.
  \item[{\rm (2)}]
  Assume that $\max\{P,Q\}>2$. Then there exists $(\mu,\nu)\in L^{r_1}({\bf R}^N)\times L^{r_2}({\bf R}^N)$
  such that problem~(P) possess no solutions for any $T>0$.
  \item[{\rm (3)}]
  If
  \begin{equation}
  \label{eq:1.1}
  \frac{q+1}{pq-1}<\frac{N}{2},
  \end{equation}
  then problem~(P) possesses a global-in-time positive solution provided that
  $(\mu,\nu)\not\equiv(0,0)$ and
  both $\|\mu\|_{L^{r_1^*,\infty}({\bf R}^N)}$ and $\|\nu\|_{L^{r_2^*,\infty}({\bf R}^N)}$ are sufficiently small,
  where
  \begin{equation}
  \label{eq:1.2}
  r_1^*:=\frac{N}{2}\frac{pq-1}{p+1},
  \qquad
  r_2^*:=\frac{N}{2}\frac{pq-1}{q+1}.
  \end{equation}
  On the other hand,
  if $(p,q)$ does not satisfy \eqref{eq:1.1}, then
  problem~(P) possesses no global-in-time positive solutions.
\end{itemize}
Unfortunately, even for the case of $D_1=D_2$,
statements~(1)--(3) are not available to the case of $(q+1)/(pq-1)\ge N/2$
and they are not enough to clarify optimal singularities of the initial functions
for the solvability of problem~(P).
On the other hand,
much less is known about the results on the solvability of problem~(P)
in the case of $D_1\not=D_2$.

Recently, in \cite{FI}
the authors of this paper studied qualitative property of the initial traces of the solutions to problem~(P)
and obtained necessary conditions on the initial data
for the existence of solutions.
More precisely, they divided problem~(P) into the following six cases:
\begin{equation*}
\begin{split}
 & {\rm (A)}\quad\frac{q+1}{pq-1}<\frac{N}{2};\\
 & {\rm (B)}\quad \frac{q+1}{pq-1}=\frac{N}{2}\quad\mbox{and}\quad p<q ;
  \qquad\qquad\,\,
	{\rm (C)}\quad \frac{q+1}{pq-1}=\frac{N}{2}\quad\mbox{and}\quad p=q ;\\
 & {\rm (D)}\quad \frac{q+1}{pq-1}>\frac{N}{2}\quad\mbox{and}\quad q>1+\frac{2}{N};
	\qquad
	{\rm (E)}\quad \frac{q+1}{pq-1}>\frac{N}{2}\quad\mbox{and}\quad q=1+\frac{2}{N};\\
 & {\rm (F)}\quad \frac{q+1}{pq-1}>\frac{N}{2}\quad\mbox{and}\quad q<1+\frac{2}{N}
\end{split}
\end{equation*}
(see Figure~1) and proved the following theorem
(see \cite[Theorems~1.1 and 1.2 ]{FI}).
\begin{theorem}
\label{Theorem:1.1}
Let $N\ge 1$, $0<p\le q$ with $pq>1$ and $T>0$.
Let $(u,v)$ be a solution to problem~{\rm (P)} in ${\bf R}^N\times[0,T)$.
Then the initial data $(\mu,\nu)$ satisfies the following:
\begin{itemize}
\item[{\rm (a)}]
  Consider case {\rm (A)}.
  Then there exists $\gamma_1>0$ such that
  $$
  \sup_{x\in{\bf R}^N}\,\mu(B(x,\sigma))\le\gamma_1 \sigma^{N-\frac{2(p+1)}{pq-1}},
  \qquad
  \sup_{x\in{\bf R}^N}\,\nu(B(x,\sigma))\le\gamma_1 \sigma^{N-\frac{2(q+1)}{pq-1}},
  $$
  for $0<\sigma\le T^\frac{1}{2};$
\item[{\rm (b)}]
  Consider case {\rm (B)}.
  Then there exists $\gamma_2>0$ such that
  $$
  \sup_{x\in{\bf R}^N}\,\int_0^\sigma \left[\frac{\mu(B(x,\tau))}{\tau^{N-\frac{2(p+1)}{pq-1}}}\right]^q\tau^{-1}\,d\tau
  +\sup_{x\in{\bf R}^N}\,\nu(B(x,\sigma))
  \le\gamma_2\biggr[\log\left(e+\frac{T^\frac{1}{2}}{\sigma}\right)\biggr]^{-\frac{1}{pq-1}}
  $$
  for $0<\sigma\le T^\frac{1}{2};$
\item[{\rm (c)}]
  Consider case {\rm (C)}.
  Then there exists $\gamma_3>0$ such that
  $$
  \sup_{x\in{\bf R}^N}\,\mu(B(x,\sigma))+\sup_{x\in{\bf R}^N}\,\nu(B(x,\sigma))
  \le\gamma_3 \left[ \log \left( e + \frac{T^\frac{1}{2}}{\sigma} \right) \right]^{-\frac{N}{2}}
  $$
  for $0<\sigma\le T^\frac{1}{2};$
\item[{\rm (d)}]
  Consider case {\rm (D)}.
  Then there exists $\gamma_4>0$ such that
  $$
  \sup_{x\in{\bf R}^N}\,\int_0^{T^\frac{1}{2}}
  \left[\frac{\mu(B(x,\tau))}{\tau^{N-\frac{N+2}{q}}}\right]^q\tau^{-1}\,d\tau
  +\sup_{x\in{\bf R}^N}\,\nu(B(x,T^\frac{1}{2})) \le \gamma_4 T^{\frac{N}{2} - \frac{q+1}{pq-1}};
  $$
\item[{\rm (e)}]
  Consider case {\rm (E)}.
  Then there exists $\gamma_5>0$ such that
  $$
  \sup_{x\in{\bf R}^N}\,\int_0^{T^\frac{1}{2}} \mu(B(x,\tau))^q\tau^{-1}\,d\tau
  +\sup_{x\in{\bf R}^N}\,\nu(B(x,T^\frac{1}{2})) \le\gamma_5 T^{\frac{N}{2} - \frac{q+1}{pq-1}};
  $$
\item[{\rm (f)}]
  Consider case {\rm (F)}.
  Then there exists $\gamma_6>0$ such that
  $$
  \sup_{x\in{\bf R}^N}\,\mu(B(x,T^\frac{1}{2})) \le\gamma_6 T^{\frac{N}{2} - \frac{p+1}{pq-1}},
  \quad
  \sup_{x\in{\bf R}^N}\,\nu(B(x,T^\frac{1}{2})) \le\gamma_6 T^{\frac{N}{2} - \frac{q+1}{pq-1}}.
  $$
\end{itemize}
Here $\gamma_1$, \dots, $\gamma_6$ are positive constants depending only on $N$, $p$, $q$, $D_1$ and $D_2$.
\end{theorem}

\begin{figure}[htbp]
	\centering
	\subfloat{
		\begin{tikzpicture}[samples = 100]
			\draw[-{Latex}] (-0.5,0) -- (5.15,0) node[right] {\small $p$};
			\draw[-{Latex}] (0,-0.5) -- (0,5.15) node[above] {\small $q$};
			\path (0,0) node[above left] {\small O};
			\fill [ opacity = 0.2 ] plot [ domain = 0.2:1 ] ({\x}, {1/\x}) -- plot [ domain = 1:5 ] ({\x}, {\x});
			\draw plot[ domain = -0.4:5 ] ({\x}, {\x});
			\path (4,4) node[below right] {\small $p=q$};
      \draw [dashed] (3/5,5/3) -- (0,5/3) node[left] {\small $1+\frac{2}{N}$};
			\draw[color = white, thick] plot[ domain = 0.2:1 ] ({\x}, {1/\x});
			\draw[dashed] plot[ domain = 0.2:5 ] ({\x}, {1/\x});
			\draw[color = white, thick] plot[ domain = 1:5/3 ] ({\x}, {(5/3) / (\x - 2/3)});
			\draw[dashed] plot[ domain = 1:5/3 ] ({\x}, {(5/3) / (\x - 2/3)});
			\draw plot[ domain = 5/3:5 ] ({\x}, {(5/3) / (\x - 2/3)});
			\path (4.2,1.2) node {\small $\dfrac{q+1}{pq-1}=\dfrac{N}{2}$};
      \path (-0.8, 3.5) node {\small $pq=1$};
      \path (-0.8, 3.2) edge [-{Latex}, bend right] (1/2.7,2.7);
      \draw [thick, color = white] (3/5,5/3) -- (5/3,5/3);
			\draw [dashed] (3/5,5/3) -- (5/3,5/3);
      \draw [dashed] (5/3,5/3) -- (5/3,0) node[below] {\small $1+\frac{2}{N}$};
			\fill [color = white] (5/3, 5/3) circle (2pt);
			\draw (5/3, 5/3) circle (2pt);
			\fill [ color = white ] (1,1) circle (2pt);
			\draw (1,1) circle (2pt);
			\path (0.85,2.8) node {\small (D)};
			\path (2.3,3.3) node[above] {\small (A)};
			\path (1.03,1.1) node[above] {\small (F)};
		\end{tikzpicture}
	}
	\quad
	\subfloat{
		\begin{tikzpicture}[samples = 100]
			\draw[-{Latex}] (-0.5,0) -- (5.15,0) node[right] {\small $p$};
			\draw[-{Latex}] (0,-0.5) -- (0,5.15) node[above] {\small $q$};
			\path (0,0) node[above left] {\small O};
			\draw plot[ domain = -0.4:5 ] ({\x}, {\x});
			\path (4,4) node[below right] {\small $p=q$};
			\draw[dashed] plot[ domain = 0.2:5 ] ({\x}, {1/\x});
			\draw plot[ domain = 1:5 ] ({\x}, {(5/3) / (\x - 2/3)});
			\path (4.2,1.2) node {\small $\dfrac{q+1}{pq-1}=\dfrac{N}{2}$};
      \path (-0.8, 3.5) node {\small $pq=1$};
      \path (-0.8, 3.2) edge [-{Latex}, bend right] (1/2.7,2.7);
      \draw [dashed] (5/3,5/3) -- (5/3,0) node[below] {\small $1+\frac{2}{N}$};
      \draw [dashed] (3/5,5/3) -- (0,5/3) node[left] {\small $1+\frac{2}{N}$};
			\fill (5/3, 5/3) circle (2pt);
			\fill [ color = white ] (1,1) circle (2pt);
			\draw (1,1) circle (2pt);
			\draw (3/5,5/3) -- (5/3,5/3);
			\path (2.8,1.9) edge [-{Latex}, bend left = 30] (1.74,5/3);
			\path (2.84,1.8) node[above] {\small (C)};
			\path (0.8,2.3) edge [-{Latex}, bend right = 15] (1,1.7);
			\path (0.8,2.2) node[above] {\small (E)};
			\path (2,3.5) edge [-{Latex}, bend left] (1.25,3);
			\path (2,3.4) node[above] {\small (B)};
		\end{tikzpicture}
	}
	\caption{}\label{figure:1}
\end{figure}
\noindent
In this paper, as a continuation of \cite{FI},
we obtain sufficient conditions on the existence of solutions to problem~(P).
Combining our main results with Theorem~\ref{Theorem:1.1},
we clarify optimal singularities of the initial functions
for the existence of solutions to problem~(P).
Our results are independent of whether $D_1=D_2$ or not.
\vspace{5pt}

Following \cite{FI}, we formulate the definition of a solution to problem~(P).
Let $G=G(x,t)$ be the fundamental solution to the heat equation in ${\bf R}^N$, that is,
\begin{equation}
\label{eq:1.3}
G(x,t)=(4\pi t)^{-\frac{N}{2}}\exp\left(-\frac{|x|^2}{4t}\right).
\end{equation}
For any Radon measure $\mu$ in $\mathbf{R}^N$, we set
\begin{equation*}
  [S(t)\mu](x) := \int_{\mathbf{R}^N} G(x-y,t)\, d\mu(y).
\end{equation*}
We also write
\begin{equation*}
  [S(t)\mu](x) = \int_{\mathbf{R}^N} G(x-y,t) \mu(y)\, dy
\end{equation*}
if $\mu$ is a nonnegative measurable function in $\mathbf{R}^N$.
%
\begin{definition}
\label{Definition:1.1}
Let $\mu$ and $\nu$ be Radon measures in ${\bf R}^N$.
Let $(u,v)$ be a pair of nonnegative measurable functions in ${\bf R}^N\times(0,T)$,
where $0<T\le\infty$.
We say that $(u,v)$ is a solution to problem~{\rm (P)} in ${\bf R}^N\times[0,T)$
if $(u,v)$ satisfies $u(x,t)<\infty$, $v(x,t)<\infty$ and
\begin{equation}
\begin{split}
\label{eq:1.4}
u(x,t)
& = [S(D_1 t)\mu](x)
+\int_0^t [S(D_1(t-s))v(s)^p](x) \,ds,\\
v(x,t)
& = [S(D_2 t)\nu](x)
+\int_0^t [S(D_2(t-s))u(s)^q](x) \,ds,
\end{split}
\end{equation}
for almost all $x\in{\bf R}^N$ and $0<t<T$.
If $(u,v)$ satisfies \eqref{eq:1.4} with $=$ replaced by $\ge$,
then we say that $(u,v)$ is a supersolution to problem~{\rm (P)}.
\end{definition}
\begin{remark}
  \label{Remark:1.1}
  {\rm (i)}
  It follows from {\rm \cite[Lemma~2.2]{IKS}}, with a slight modification, that
  problem~{\rm (P)} possesses a solution in $\mathbf{R}^N \times [0,T)$
  if and only if problem~{\rm (P)} possesses a supersolution in $\mathbf{R}^N \times [0,T)$.
  \vspace{3pt}
  \newline
  {\rm (ii)} Let $T>0$ and let $(u,v)$ be a solution to problem~{\rm (P)} in ${\bf R}^N\times[0,T)$.
  Let $\alpha>0$. Set
  $$
  \hat{u}(x,t):=T^{\frac{p+1}{pq-1}}u(\alpha T^{1/2}x,Tt),
  \qquad
  \hat{v}(x,t):=T^{\frac{q+1}{pq-1}}v(\alpha T^{1/2}x,Tt),
  $$
  for $x\in {\bf R}^N$ and $t\in(0,1)$.
  Then $(\hat{u},\hat{v})$ satisfies
  \begin{equation*}
  \left\{
  \begin{array}{ll}
  \partial_t \hat{u}=D_1\alpha^{-2}\Delta\hat{u}+\hat{v}^p & \quad\mbox{in}\quad{\bf R}^N\times(0,1),\vspace{3pt}\\
  \partial_t \hat{v}=D_2\alpha^{-2}\Delta\hat{v}+\hat{u}^q & \quad\mbox{in}\quad{\bf R}^N\times(0,1),\vspace{3pt}\\
  \hat{u},\hat{v}\ge 0 & \quad\mbox{in}\quad{\bf R}^N\times(0,T),\vspace{3pt}\\
  (\hat{u}(\cdot,0),\hat{v}(\cdot,0))=(\hat{\mu},\hat{\nu}) & \quad\mbox{in}\quad{\bf R}^N.
  \end{array}
  \right.
  \end{equation*}
  Here $\hat{\mu}$ and $\hat{\nu}$ are Radon measure in ${\bf R}^N$ such that
  $$
  \hat{\mu}(K)=\alpha^{-N}T^{\frac{p+1}{pq-1}-\frac{N}{2}}\mu(\alpha T^{\frac{1}{2}}K),
  \quad
  \hat{\nu}(K)=\alpha^{-N}T^{\frac{q+1}{pq-1}-\frac{N}{2}}\nu(\alpha T^{\frac{1}{2}}K),
  $$
  for Borel sets $K$ in ${\bf R}^N$.
\end{remark}

We are ready to state one of the main results of this paper.
By Theorem~\ref{Theorem:1.2} we clarify optimal singularities of the initial functions
for the solvability of problem~(P).
\begin{Theorem}
  \label{Theorem:1.2}
  Let $N\ge 1$ and $0<p\le q$ with $pq>1$.
  \begin{itemize}
    \item[{\rm (a)}]
    Consider case {\rm (A)}. Let
    \begin{align*}
      & \mu(x)=c_{a,1} |x|^{-\frac{2(p+1)}{pq-1}}\chi_{B(0,1)}(x)\quad\mbox{in}\quad{\bf R}^N,\\
      & \nu(x)=c_{a,2} |x|^{-\frac{2(q+1)}{pq-1}}\chi_{B(0,1)}(x)\quad\mbox{in}\quad{\bf R}^N,
    \end{align*}
    where $c_{a,1}$, $c_{a,2}>0$.
    Then
    problem~{\rm (P)} possesses no positive local-in-time solutions
    if either $c_{a,1}$ or $c_{a,2}$ is sufficiently large.
    On the other hand,
    problem~{\rm (P)} possesses a global-in-time solution
    if both of $c_{a,1}$ and $c_{a,2}$ are sufficiently small.
    \item[{\rm (b)}]
    Consider case {\rm (B)}. Let
    \begin{align*}
      & \mu(x)=c_{b,1} |x|^{-\frac{2(p+1)}{pq-1}}\left|\log\frac{|x|}{2}\right|^{-\frac{p}{pq-1}}\chi_{B(0,1)}(x)
      \quad\mbox{in}\quad{\bf R}^N,\\
      & \nu(x)=c_{b,2} |x|^{-N}\left|\log\frac{|x|}{2}\right|^{-\frac{1}{pq-1}-1}\chi_{B(0,1)}(x)
     \quad\,\,\,\,\mbox{in}\quad{\bf R}^N,
    \end{align*}
    where $c_{b,1}$, $c_{b,2}>0$.
    Then problem~{\rm (P)} possesses no positive local-in-time solutions
    if either $c_{b,1}$ or $c_{b,2}$ is sufficiently large.
    On the other hand,
    problem~{\rm (P)} possesses a local-in-time solution
    if both of $c_{b,1}$ and $c_{b,2}$ are sufficiently small.
    \item[{\rm (c)}]
    Consider case {\rm (C)}. Let
    \begin{align*}
      & \mu(x)=c_{c,1}|x|^{-N}\left|\log\frac{|x|}{2}\right|^{-\frac{N}{2}-1}\chi_{B(0,1)}(x)\quad\mbox{in}\quad{\bf R}^N,\\
      & \nu(x)=c_{c,2}|x|^{-N}\left|\log\frac{|x|}{2}\right|^{-\frac{N}{2}-1}\chi_{B(0,1)}(x)\quad\mbox{in}\quad{\bf R}^N,
    \end{align*}
    where $c_{c,1}$, $c_{c,2}>0$.
    Then problem~{\rm (P)} possesses no positive local-in-time solutions
    if either $c_{c,1}$ or $c_{c,2}$ is sufficiently large.
    On the other hand,
    problem~{\rm (P)} possesses a local-in-time solution
    if both of $c_{c,1}$ and $c_{c,2}$ are sufficiently small.
    \item[{\rm (d)}]
    Consider case {\rm (D)}.
    Let
    $$
    \mu(x)=|x|^{-\frac{N+2}{q}}h_1(|x|)\chi_{B(0,1)}(x)\quad\mbox{in}\quad{\bf R}^N,
    $$
    where $h_1$ is a positive increasing function in $(0,1]$
    such that $h_1(1)<\infty$
    and $r^{-\epsilon}h_1(r)$ is decreasing in $r$ for some $\epsilon>0$.
    Let $\nu$ be a Radon measure in ${\bf R}^N$.
    Then problem~{\rm (P)} possesses no positive local-in-time solution if either
    $$
    \int_0^1 h_1(\tau)^q\tau^{-1}\,d\tau=\infty
    \quad\mbox{or}\,\,\,
    \sup_{x\in{\bf R}^N}\nu(B(x,1))=\infty.
    $$
    On the other hand,
    problem~{\rm (P)} possesses a local-in-time solution if
    $$
    \int_0^1 h_1(\tau)^q\tau^{-1}\,d\tau<\infty
    \quad\mbox{and}\,\,\,
    \sup_{x\in{\bf R}^N}\nu(B(x,1))<\infty.
    $$
    \item[{\rm (e)}]
    Consider case {\rm (E)}.
    Let
    $$
    \mu(x)=|x|^{-N}h_2(|x|)\chi_{B(0,1)}(x)\quad\mbox{in}\quad{\bf R}^N,
    $$
    where $h_2$ is a positive increasing function in $(0,1]$ satisfying $h_2(1)<\infty$.
    Let $\nu$ be a Radon measure in ${\bf R}^N$.
    Then problem~{\rm (P)} possesses no positive local-in-time solutions if either
    $$
    \int_0^1\left[\int_0^r h_2(\tau)\tau^{-1}\,d\tau\right]^qr^{-1}\,dr=\infty
    \quad\mbox{or}\,\,\,
    \sup_{x\in{\bf R}^N}\nu(B(x,1))=\infty.
    $$
    On the other hand,
    problem~{\rm (P)} possesses a local-in-time solution if
    $$
    \int_0^1\left[\int_0^r h_2(\tau)\tau^{-1}\,d\tau\right]^qr^{-1}\,dr<\infty
    \quad\mbox{and}\,\,\,
    \sup_{x\in{\bf R}^N}\nu(B(x,1))<\infty.
    $$
    \item[{\rm (f)}]
    Consider case {\rm (F)}.
    Let $\mu$ and $\nu$ be Radon measures in ${\bf R}^N$.
    Then problem~{\rm (P)} possesses no positive local-in-time solutions if either
    $$
    \sup_{x\in{\bf R}^N}\mu(B(x,1))=\infty
    \quad\mbox{or}\,\,\,
    \sup_{x\in{\bf R}^N}\nu(B(x,1))=\infty.
    $$
    On ther other hand, problem~{\rm (P)} possesses a local-in-time solution if
    $$
    \sup_{x\in{\bf R}^N}\mu(B(x,1))<\infty
    \quad\mbox{and}\,\,\,
    \sup_{x\in{\bf R}^N}\nu(B(x,1))<\infty.
    $$
  \end{itemize}
\end{Theorem}
The statements of the nonexistence of local-in-time solutions in Theorem~\ref{Theorem:1.2}
follow from Theorem~\ref{Theorem:1.1}. (See also \cite[Corollary~1.2]{FI}.)
The other statements in Theorem~\ref{Theorem:1.2} follow from the results in Sections~3--5,
which are the main ingredients of this paper and which are concerned with sufficient conditions
on the solvability of problem~(P) with singular initial functions.

Optimal singularities of the initial functions for the solvability of problem~(P)
depend on $(p,q)$ and so do our sufficient conditions on the existence of solutions.
Problem~(P) in cases (A), (C) and (F) can be regarded as generalizations of
the Cauchy problem for
$$
\partial_t u=\Delta u+u^p\quad\mbox{in}\quad{\bf R}^N\times(0,T)
$$
in the cases $p>1+2/N$, $p=1+2/N$ and $p>1+2/N$, respectively.
We construct a supersolution in caces (A), (C) and (F)
by applying similar arguments in \cite{HI01} and \cite{IKS}.
Then, thanks to Remark~\ref{Remark:1.1}~(i),
we give sufficient conditions for the existence of solutions for problem~\eqref{eq:P}
(see Theorems~\ref{Theorem:3.1}, \ref{Theorem:3.2} and \ref{Theorem:3.3} in Section~3).
Cases (B), (D) and (E) are specific to the parabolic system
and the construction of supersolutions
requires delicate decay estimates of $S(t)\mu$ (see lemmas in Section~2)
and more complicate arguments than those in cases (A), (C) and (F)
(see Theorems~\ref{Theorem:4.1} and \ref{Theorem:5.1}).

The rest of this paper is organized as follows.
In Section~2 we recall some properties on uniformly local weak $L^r$ spaces
and obtain some decay estimates of $S(t)\varphi$.
Section~3 is devoted to cases (A), (C) and (F).
In Sections~4 and 5 we consider case (B) and cases (D) and (E), respectively.
In Section~6 we apply the results in Sections~3, 4 and 5 to problem~\eqref{eq:P}
and prove Theorem~\ref{Theorem:1.2}.
This shows the validity of our sufficient conditions given in Sections~3, 4 and 5
for the existence of solutions to problem~(P).
\section{Preliminaries}
In this section we introduce some notations and
prove several lemmas on $S(t)\varphi$.
In what follows, we denote by $C$ a generic constant depending only on
$p$, $q$, $D_1$, $D_2$ and $N$, which may change line by line.

We introduce some notations.
For any measurable set $\Omega$ in ${\bf R}^N$ and $1\le r\le\infty$,
$L^r(\Omega)$ denotes the usual Lebesgue space equipped with the norm $\|\cdot\|_{L^r(\Omega)}$.
In the case of $\Omega={\bf R}^N$,
we write $L^r=L^r({\bf R}^N)$ and $\|\cdot\|_{L^r}=\|\cdot\|_{L^r({\bf R}^N)}$ for simplicity.
Let $L^r_{\rm loc}$ be the local $L^r$ space in ${\bf R}^N$.
We introduce the uniformly local Lebesuge space $L_{{\rm uloc}}^r$ by
$$
L_{{\rm uloc}}^r:=
\left\{f\,:\,\mbox{$f\in L^r_{{\rm loc}}$ such that $|||f|||_r<\infty$}\right\},
$$
where
$$
|||f|||_r:=\sup_{x\in{\bf R}^N}\|f\|_{L^r(B(x,1))}.
$$
Notice that $|||f||_r<\infty$ if and only if
$$
|||f|||_{r,\rho}:=\sup_{x\in{\bf R}^N}\|f\|_{L^r(B(x,\rho))}<\infty,
\quad \rho>0.
$$

We recall some properties of $S(t)\varphi$.
Let $\varphi$ be a nonnegative measurable function in $\mathbf{R}^N$
and $\Phi$ an increasing convex function in $[0,\infty)$.
It follows from the Jensen inequality that
\begin{gather}
  \label{eq:2.1}
  [S(t)\varphi](x) \le \Phi^{-1}\left( [S(t)\Phi(\varphi)](x) \right),
  \quad x\in \mathbf{R}^N,\,\,t>0.
\end{gather}
The following inequalities hold:
\begin{align}
 \label{eq:2.2}
  &  \|S(t)\varphi\|_{L^\ell}
  \le Ct^{-\frac{N}{2}\left(\frac{1}{r}-\frac{1}{\ell}\right)}\|\varphi\|_{L^r},
  \quad \varphi\in L^r, \\
\label{eq:2.3}
  & |||S(t)\varphi|||_{\ell,\rho}
    \le C
    \left[ t^{-\frac{N}{2}\left(\frac{1}{r}-\frac{1}{\ell}\right)} + \rho^{-N \left(\frac{1}{r}-\frac{1}{\ell}\right)} \right]
    |||\varphi|||_{r,\rho},\quad\varphi\in L^r_{{\rm uloc}},
\end{align}
for $t>0$, $\rho>0$ and $1\le r\le\ell\le\infty$
(see e.g., \cite[(1.8)]{MT} for \eqref{eq:2.3}).
In particular,
\begin{equation}
\label{eq:2.4}
||S(t)\varphi|||_\ell
\le Ct^{-\frac{N}{2}\left(\frac{1}{r}-\frac{1}{\ell}\right)}|||\varphi|||_r,\quad\varphi\in L^r_{{\rm uloc}}
\end{equation}
for $0<t\le 1$ and $1\le r\le\ell\le\infty$.
Furthermore, we have:
\begin{lemma}
\label{Lemma:2.1}
There exists $C>0$ such that
\begin{equation}
\label{eq:2.5}
\|S(t)\mu\|_\infty\le Ct^{-\frac{N}{2}}
\sup_{x\in{\bf R}^N}\mu(B(x,\sqrt{t})),\quad t>0,
\end{equation}
for nonnegative Radon measures $\mu$ in ${\bf R}^N$. 
\end{lemma}
See also \cite[Lemma~2.1]{HI01}.
\vspace{5pt}

We next obtain estimates of $S(t)\mu$ in an annular domain.
\begin{Lemma}
  \label{Lemma:2.2}
  Let $\mu$ be a nonnegative measurable function in $\mathbf{R}^N$.
  Assume that there exist $a\in (0,N]$
  and a nonnegative increasing function $f$ on $(0,1]$
  such that $f(1)<\infty$ and
  \begin{equation}
    \notag
    \mu(x) \le |x|^{-a} f(|x|) \chi_{B(0,1)}(x) \quad\mbox{in}\,\,\, \mathbf{R}^N.
  \end{equation}
  Then, for any $r_* > N/a$,
  there exists $C>0$ such that
  \begin{equation}
    \label{eq:2.6}
    \| S(t)\mu \|_{L^{r_*}(B(0,1)\setminus B(0,\sqrt{t}))}
    \le C t^{-\frac{N}{2}\left( \frac{a}{N} - \frac{1}{r_*} \right)} g(t)
  \end{equation}
  for sufficiently small $t>0$,
  where $g$ is a function defined by
  \begin{equation}
    \label{eq:2.7}
    g(t) = \left\{
    \begin{array}{ll}
      f(t^\frac{1}{6}) + t^\frac{ar_*-N}{4r_*}
      &\mbox{if}\quad a<N, \\[5pt]
      f(t^\frac{1}{6}) + t^\frac{Nr_*-N}{4r_*}
      + \displaystyle\int_0^{\sqrt{t}} \tau^{-1}f(\tau)\, d\tau
      &\mbox{if}\quad a=N.
    \end{array}
    \right.
  \end{equation}
\end{Lemma}
\textbf{Proof.}
Let $t\in (0,1)$ and $|x|\ge \sqrt{t}$.
By \eqref{eq:1.3} we have
\begin{align*}
  [S(t)\mu](x)
  & =
  (4\pi t)^{-\frac{N}{2}}\left[
  \int_{\{|y-x|\le |x|/2\}} e^{-\frac{|x-y|^2}{4t}}\mu(y)\, dy
  + \int_{\{|y-x|\ge |x|/2\}} e^{-\frac{|x-y|^2}{4t}}\mu(y)\, dy
  \right]
  \\
  & =:
  I(x,t) + J(x,t).
\end{align*}
Since $|y|\ge |x|/2$ and $|y|\le 3|x|/2$ if
$|y-x|\le |x|/2$,
we have
\begin{align*}
  I(x,t)
  & \le
  C(4\pi t)^{-\frac{N}{2}} \int_{\{|y-x|\le |x|/2\}} e^{-\frac{|y-x|^2}{4t}} \, dy
  \cdot |x|^{-a} f\left(\frac{3|x|}{2}\right)
  \\
  & \le
  C |x|^{-a} f\left(\frac{3|x|}{2}\right)
  \int_{\{|\eta|\le |x|/2\sqrt{t}\}} e^{-\frac{|\eta|^2}{4}} \, d\eta
  \le
  C |x|^{-a} f\left(\frac{3|x|}{2}\right).
\end{align*}
Then we obtain
\begin{equation}
  \label{eq:2.8}
  \begin{aligned}
    \int_{\{\sqrt{t}\le |x|\le t^\frac{1}{4}\}} |I(x,t)|^{r_*} \, dx
    & \le
    C\int_{\{1\le |x|\le t^{-\frac{1}{4}}\}} |x|^{-ar_*} f\left(\frac{3\sqrt{t} |x|}{2}\right)^{r_*} \, dx
    \cdot t^{-\frac{ar_*}{2} + \frac{N}{2}}
    \\
    & \le
    C t^{-\frac{ar_*}{2} + \frac{N}{2}} f\left(\frac{3t^\frac{1}{4}}{2}\right)^{r_*}
    \int_{\{|x|\ge 1\}} |x|^{-ar_*}\, dx
    \\
    & \le
    C t^{-\frac{ar_*}{2} + \frac{N}{2}} f(t^\frac{1}{6})^{r_*}
    \int_{\{|x|\ge 1\}} |x|^{-ar_*}\, dx
  \end{aligned}
\end{equation}
for sufficiently small $t>0$.
Note that the last integral converges due to $ar_*>N$.
Similarly, we obtain
\begin{equation}
  \label{eq:2.9}
  \begin{aligned}
    \int_{\{t^\frac{1}{4}\le |x|\le 1\}} | I(x,t) |^{r_*} \, dx
    & \le
    C f(1) \int_{\{1\le |x|\le t^{-\frac{1}{4}}\}}
    |x|^{-ar*} \, dx \cdot t^{-\frac{ar_*}{4}+\frac{N}{4}}
    \\
    & \le
    C f(1) t^{-\frac{ar_*}{2} + \frac{N}{2}} \cdot t^{\frac{ar_*}{4} - \frac{N}{4}} \int_{\{|x|\ge 1\}} |x|^{-ar_*}\, dx.
  \end{aligned}
\end{equation}
Therefore, by \eqref{eq:2.8} and \eqref{eq:2.9} we have
\begin{equation}
  \label{eq:2.10}
  \left[ \int_{\{\sqrt{t} \le |x|\le 1\}}
  | I(x,t) |^{r_*} \, dx \right]^{r_*}
  \le C t^{-\frac{N}{2} \left( \frac{a}{N} - \frac{1}{r_*} \right)}
  \left[ f(t^\frac{1}{6}) + t^{\frac{ar_* - N}{4r_*}} \right]
\end{equation}
for sufficiently small $t>0$.

On the other hand,
since $|y-x| \ge |x|/2$ if $|x|\ge \sqrt{t}$ and $|y| \le \sqrt{t}/2$,
we see that
\begin{align*}
   & J(x,t)
  =
  (4\pi t)^{-\frac{N}{2}} \left[
  \int_{\{|y|\le \sqrt{t}/2\}} + \int_{ \{ |y|\ge \frac{\sqrt{t}}{2} \} \cap \{ |y-x|\ge |x|/2\} }
  \right]
  e^{-\frac{|y-x|^2}{4t}} \mu(y)\, dy
  \\
  & \le
  (4\pi t)^{-\frac{N}{2}} \int_{\{|y|\le \sqrt{t}/2\}} e^{-\frac{|x|^2}{16t}} \mu(y)\, dy
  + (4\pi t)^{-\frac{N}{2}} \int_{ \{ |y|\ge \sqrt{t}/2 \} \cap \{ |y-x|\ge |x|/2\} }
  e^{-\frac{|y-x|^2}{4t}} \mu(y)\, dy
  \\
  & =:
  J_1(x,t) + J_2(x,t).
\end{align*}
We now divide the proof into two cases $a<N$ and $a=N$.

\smallskip

\noindent
\underline{Case $a<N$} :
Since $|x|\ge 2|y|$ if $|x|\ge \sqrt{t}$ and $|y| \le \sqrt{t}/2$, we have
\begin{align*}
  J_1(x,t)
  & \le (4\pi t)^{-\frac{N}{2}}
  \int_{\{|y|\le \sqrt{t}/2\}} |y|^{-a} \, dy
  \cdot e^{-\frac{|x|^2}{16t}}f\left( \frac{\sqrt{t}}{2} \right)
  \\
  & \le C t^{-\frac{a}{2}}
  e^{-\frac{|x|^2}{16t}} f(\sqrt{t})
  \int_{\{|y|\le 1/2\}} |y|^{-a} \, dy.
\end{align*}
The last integral converges due to $a<N$.
Then we have
\begin{equation}
  \label{eq:2.11}
  \begin{aligned}
    \int_{\{\sqrt{t}\le |x|\le 1\}} |J_1(x,t)|^{r_*} \, dx
    & \le
    C t^{-\frac{ar_*}{2}} f(\sqrt{t})^{r_*} \int_{\{|x|\le 1\}} e^{-\frac{r_* |x|^2}{16t}} \, dx
    \\
    & \le
    C t^{-\frac{ar_*}{2} + \frac{N}{2}} f(\sqrt{t})^{r_*}.
  \end{aligned}
\end{equation}
On the other hand, since
\begin{equation*}
  \mu(y) \le
  \left\{
  \begin{array}{cll}
    C t^{-\frac{a}{2}} f( t^\frac{1}{4} ) &\mbox{if} &\frac{\sqrt{t}}{2}\le |y|\le t^\frac{1}{4}, \\[5pt]
    C t^{-\frac{a}{4}} f(1) &\mbox{if} &|y|\ge t^\frac{1}{4},
  \end{array}
  \right.
\end{equation*}
we have
\begin{align*}
  J_2(x,t)
  & \le
  C t^{-\frac{N}{2}} \int_{\{|y-x|\ge \frac{|x|}{2}\}}
  e^{-\frac{|y-x|^2}{4t}} \, dy \cdot
  \left[ t^{-\frac{a}{2}} f(t^\frac{1}{4}) + t^{-\frac{a}{4}}f(1) \right]
  \\
  & \le
  C t^{-\frac{N}{2}} \int_{\mathbf{R}^N} e^{-\frac{|y-x|^2}{8t}} \, dy
  \cdot \left[ t^{-\frac{a}{2}} f(t^\frac{1}{4}) + t^{-\frac{a}{4}}f(1) \right]
  \cdot e^{-\frac{|x|^2}{32t}}\\
   & \le
  C\left[ t^{-\frac{a}{2}} f(t^\frac{1}{4}) + t^{-\frac{a}{4}}f(1) \right]
  e^{-\frac{|x|^2}{32t}}.
\end{align*}
Then we see that
\begin{equation}
  \label{eq:2.12}
  \begin{split}
    & \int_{\{\sqrt{t}\le |x|\le 1\}} |J_2(x,t)|^{r_*}\, dx
    \\
    & \le
    C \left[ t^{-\frac{ar_*}{2}} f(\sqrt{t})^{r_*} + t^{-\frac{ar_*}{4}}f(1)^{r_*} \right]
    \int_{\{\sqrt{t}\le |x|\le 1\}} e^{-\frac{r_* |x|^2}{32t}} \, dx
    \\
    & \le
    C \left[ t^{-\frac{ar_*}{2} + \frac{N}{2}} f(\sqrt{t})^{r_*} + t^{-\frac{ar_*}{4} + \frac{N}{2}}f(1)^{r_*} \right]
  \end{split}
\end{equation}
for sufficiently small $t>0$.
Therefore, by \eqref{eq:2.11} and \eqref{eq:2.12} we obtain
\begin{equation}
  \label{eq:2.13}
  \left[ \int_{\{\sqrt{t}\le |x|\le 1\}} |J(x,t)|^{r_*} \right]^\frac{1}{r_*}
  \le Ct^{-\frac{N}{2}\left( \frac{a}{N} - \frac{1}{r_*} \right)}
  \left[ f(\sqrt{t}) + t^{\frac{a}{4}} \right]
\end{equation}
for sufficiently small $t>0$.

\smallskip

\noindent
\underline{Case $a=N$} :
Since $|x|\ge 2|y|$ if $|x|\ge \sqrt{t}$ and $|y| \le \sqrt{t}/2$, we have
\begin{align*}
  J_1(x,t)
  & \le
  (4\pi t)^{-\frac{N}{2}} \int_{\{|y|\le \sqrt{t}/2\}} |y|^{-N} f(|y|) \, dy \cdot e^{-\frac{|x|^2}{16t}}
  \\
  & \le
  C t^{-\frac{N}{2}} e^{-\frac{|x|^2}{16t}} \int_0^{\frac{\sqrt{t}}{2}} \tau^{-1} f(\tau) \, d\tau
  \le C t^{-\frac{N}{2}} e^{-\frac{|x|^2}{16t}} \int_0^{\sqrt{t}} \tau^{-1} f(\tau) \, d\tau.
\end{align*}
Setting
\begin{equation*}
  \tilde{f}(s)
  := \int_0^{s} \tau^{-1} f(\tau) \, d\tau,
  \qquad 0<s\le 1,
\end{equation*}
as in \eqref{eq:2.11}, we obtain
\begin{equation}
  \label{eq:2.14}
  \int_{\{\sqrt{t}\le |x|\le 1\}} |J_1(x,t)|^{r_*} \, dx
  \le
  C t^{-\frac{Nr_*}{2} + \frac{N}{2}} \tilde{f}(\sqrt{t})^{r_*}.
\end{equation}
On the other hand, we can derive the same calculation for $J_2$ as in the case $a<N$,
and by \eqref{eq:2.12} and \eqref{eq:2.14} we obtain
\begin{equation}
  \label{eq:2.15}
  \left[ \int_{\{\sqrt{t}\le |x|\le 1\}} |J(x,t)|^{r_*}\, dx \right]^\frac{1}{r_*}
  \le C t^{-\frac{N}{2}\left( 1 - \frac{1}{r_*} \right)}
  \left[ f(\sqrt{t}) + \tilde{f}(\sqrt{t}) + t^{\frac{a}{4}} \right]
\end{equation}
for sufficiently small $t>0$.

Thus, since $f$ is an increasing function in $(0,1]$ and $\sqrt{t}<t^{1/6}$ for $0<t<1$,
by \eqref{eq:2.10}, \eqref{eq:2.13} and \eqref{eq:2.15}
we obtain inequality~\eqref{eq:2.6}.
Thus Lemma~\ref{Lemma:2.2} follows.
$\Box$
\vspace{5pt}

At the end of this section we prove a lemma,
which gives an upper bound of some integrals including logarithmic functions.
\begin{lemma}
	\label{Lemma:2.3}
	Let $a>-1$ and $b\in \mathbf{R}$.
	Then there exists $C>0$ such that
  \begin{equation*}
    \int_0^t s^a \left|\log \frac{s}{2}\right|^b \, ds
  	\le Ct^{a+1} \left|\log \frac{t}{2}\right|^b,\quad 0<t<1.
  \end{equation*}
\end{lemma}
\textbf{Proof.}
Set
\begin{equation*}
  I(b,t):= \int_0^t s^a \left|\log \frac{s}{2}\right|^b \, ds
\end{equation*}
for $0<t<1$.
If $b<0$, then we have
\begin{equation}
  \label{eq:2.16}
  \int_0^t s^a \left|\log \frac{s}{2}\right|^b \, ds
  \le \left|\log \frac{t}{2}\right|^b \int_0^t s^a \, ds
  = \frac{1}{a+1}t^{a+1}\left|\log \frac{t}{2}\right|^b
\end{equation}
for $0<t<1$.
For $b\ge 0$, by integration by parts we have
\begin{equation}
  \label{eq:2.17}
	I(b,t)
  = \int_0^t \left( \frac{1}{a+1} s^{a+1} \right)' \left|\log \frac{s}{2}\right|^b \, ds
  = \frac{t^{a+1}}{a+1}\left|\log \frac{t}{2}\right|^b
  + \frac{b}{a+1} I(b-1)
\end{equation}
for $0<t<1$.
Repeating the above argument, we see that
\begin{equation}
  \notag
  I(b,t) \le C t^{a+1}\left|\log \frac{t}{2}\right|^b + C I(b-[b]-1,t)
\end{equation}
for $0<t<1$,
where $[b]$ denote the greatest integer less than or equal to $b$.
Since $b-[b]-1<0$,
by \eqref{eq:2.16} and \eqref{eq:2.17} we obtain the desired inequality.
Thus Lemma~\ref{Lemma:2.3} follows.
$\Box$
\section{Cases (A), (C) and (F)}
In this section we focus on cases (A), (C) and (F) and obtain sufficient conditions
on the existence of solutions to problem~(P).
\begin{theorem}
\label{Theorem:3.1}
Let $N\ge 1$ and $0<p\le q$ with $pq>1$ be in case {\rm (A)}.
Let
\begin{equation}
\label{eq:3.1}
1<\alpha<\frac{pq+q}{q+1}.
\end{equation}
Then there exists $\gamma>0$ such that,
if $\mu$ and $\nu$ are nonnegative measurable functions in ${\bf R}^N$ and satisfy
\begin{equation}
\label{eq:3.2}
\big\|S(t)\mu^{\frac{\alpha(q+1)}{p+1}}\big\|_\infty
+ \big\|S(t)\nu^\alpha\big\|_\infty
\le \gamma t^{-\frac{q+1}{pq-1}\alpha},\quad 0<t<1,
\end{equation}
then problem~{\rm (P)} possesses  a solution in ${\bf R}^N\times[0,1)$.
\end{theorem}
\begin{theorem}
\label{Theorem:3.2}
Let $N\ge 1$ and $0<p\le q$ with $pq>1$ be in case {\rm (C)}.
Let $\beta>0$ and set $\Phi(\tau):=\tau[\log(e+\tau)]^\beta$ for $\tau\ge 0$.
Then there exists $\gamma>0$ such that,
if $\mu$ and $\nu$ are nonnegative measurable functions in ${\bf R}^N$ and satisfy
\begin{equation}
\label{eq:3.3}
\|S(t)\Phi(\mu)\|_\infty
+\|S(t)\Phi(\nu)\|_\infty
\le\gamma t^{-\frac{N}{2}}\left|\log\frac{t}{2}\right|^{-\frac{N}{2}+\beta},
\quad 0<t<1,
\end{equation}
then problem~{\rm (P)} possesses a solution in ${\bf R}^N\times[0,1)$.
\end{theorem}
\begin{theorem}
\label{Theorem:3.3}
Let $N\ge 1$ and $0<p\le q$ with $pq>1$ be in case {\rm (F)}.
Then there exists $\gamma>0$ such that,
if $\mu$ and $\nu$ are Radon measures in ${\bf R}^N$
and satisfy
\begin{equation}
\label{eq:3.4}
\|S(t)\mu\|_\infty+\|S(t)\nu\|_\infty
\le\gamma t^{-\frac{N}{2}},\quad 0<t<1,
\end{equation}
then problem~{\rm (P)} possesses a solution in ${\bf R}^N\times[0,1)$.
\end{theorem}
Let $D := \min\{D_1,D_2\}$ and $D' := \max\{D_1,D_2\}$.
Due to Remark~\ref{Remark:1.1}~(ii),
it suffices to consider the case where
\begin{equation}
\label{eq:3.5}
T=1,\qquad 0<D\le D'=\max\{D_1,D_2\}=1.
\end{equation}
We construct supersolutions to problem~(P) and
prove Theorems~\ref{Theorem:3.1}, \ref{Theorem:3.2} and \ref{Theorem:3.3}.
It follows that
\begin{equation}
  \label{eq:3.6}
  G(x,D_i t) = (4\pi D_i t)^{-\frac{N}{2}}
  \exp\left( -\frac{|x|^2}{4D_i t} \right)
  \le D^{-\frac{N}{2}} G(x,D' t)
  =D^{-\frac{N}{2}} G(x,t)
\end{equation}
for $x\in \mathbf{R}^N$ and $t>0$, where $i\in \{1,2\}$.
Let $(\tilde{u},\tilde{v})$ be a solution
to the Cauchy problem
\begin{equation}
  \tag{P'}
  \left\{
  \begin{array}{lll}
    \partial_t u =\Delta u +D^{-\frac{N}{2}}v^p
    & \mbox{in} & \mathbf{R}^N \times (0,1), \\[3pt]
    \partial_t v =\Delta v + D^{-\frac{N}{2}}u^q
    & \mbox{in} & \mathbf{R}^N \times (0,1), \\[3pt]
    u,v \ge 0 & \mbox{in} & \mathbf{R}^N \times (0,1), \\[3pt]
    (u(0),v(0)) =(\mu_D, \nu_D)
    & \mbox{in} & \mathbf{R}^N,
  \end{array}
  \right.
\end{equation}
where $(\mu_D,\nu_D):=D^{-\frac{N}{2}}(\mu,\nu)$.
By Definition~\ref{Definition:1.1}, \eqref{eq:3.5} and  \eqref{eq:3.6} we see that
\begin{equation*}
\begin{split}
\tilde{u}(x,t) & =D^{-\frac{N}{2}}[S(t)\mu](x)+D^{-\frac{N}{2}}\int_0^t [S(t-s)\tilde{v}(s)^p](x)\,ds\\
 & \ge [S(D_1t)\tilde{u}(0)](x)+\int_0^t [S(D_1(t-s))\tilde{v}(s)^p](x)\,ds
\end{split}
\end{equation*}
for almost all ${\bf R}^N\times(0,1)$. Similarly, we have
$$
\tilde{v}(x,t)\ge [S(D_2t)\tilde{v}(0)](x)+\int_0^tS(D_2(t-s))\tilde{u}(s)^q](x)\,ds
$$
for almost all ${\bf R}^N\times(0,1)$.
This implies that $(\tilde{u},\tilde{v})$ is a supersolution to problem~\eqref{eq:P}.
By Remark~\ref{Remark:1.1}~(i)
we see that problem~(P) possesses a solution in ${\bf R}^N\times[0,1)$
if there exists a solution to problem~(P') in ${\bf R}^N\times[0,1)$.
\vspace{5pt}

\noindent
{\bf Proof of Theorem~\ref{Theorem:3.1}.}
It suffices to construct a supersolution to problem~(P') in ${\bf R}^N\times[0,1)$.
Set
\begin{equation}
\label{eq:3.7}
\begin{split}
  & w(x,t):=\left[S(t)\mu_D^{\frac{\alpha(q+1)}{p+1}}\right](x)+[S(t)\nu_D^\alpha](x),\\
  & \overline{u}(x,t):=2w(x,t)^{\frac{p+1}{\alpha(q+1)}},
 \quad\overline{v}(x,t):=2w(x,t)^{\frac{1}{\alpha}}.
\end{split}
\end{equation}
Then
\begin{equation}
\label{eq:3.8}
0\le \mu_D(x)\le\overline{u}(x,0),\quad 0\le \nu_D(x)\le\overline{v}(x,0),\qquad x\in{\bf R}^N.
\end{equation}
Furthermore, it follows from \eqref{eq:3.2} that
\begin{equation}
\label{eq:3.9}
\|w(t)\|_\infty\le C\gamma t^{-\frac{q+1}{pq-1}\alpha},\qquad 0<t<1.
\end{equation}
By the Jensen inequality (see \eqref{eq:2.1}), \eqref{eq:3.1} and \eqref{eq:3.7} we have
\begin{equation}
\label{eq:3.10}
\begin{split}
 & [S(t)\nu_D](x)+D^{-\frac{N}{2}}\int_0^t [S(t-s)\overline{u}(s)^q](x)\,ds\\
 & \le w(x,t)^{\frac{1}{\alpha}}+C\int_0^t\left[S(t-s)w(s)^{\frac{pq+q}{\alpha(q+1)}}\right](x)\,ds\\
 & \le w(x,t)^{\frac{1}{\alpha}}+C\int_0^t \|w(s)\|_\infty^{\frac{pq+q}{\alpha(q+1)}-1}[S(t-s)w(s)](x)\,ds\\
 & \le w(x,t)^{\frac{1}{\alpha}}+Cw(x,t)\int_0^t \|w(s)\|_\infty^{\frac{pq+q}{\alpha(q+1)}-1}\,ds\\
 & \le w(x,t)^{\frac{1}{\alpha}}+C\|w(t)\|_\infty^{1-\frac{1}{\alpha}}w(x,t)^{\frac{1}{\alpha}}
 \int_0^t \|w(s)\|_\infty^{\frac{pq+q}{\alpha(q+1)}-1}\,ds
\end{split}
\end{equation}
for $(x,t)\in{\bf R}^N\times(0,1)$.
On the other hand, it follows from $\alpha>1$ that
\begin{equation*}
  -\frac{q+1}{pq-1}\alpha \left( \frac{pq+q}{\alpha(q+1)} - 1 \right)
  = \frac{-pq-q+\alpha (q+1)}{pq-1}
  > -1.
\end{equation*}
Then, by \eqref{eq:3.9} and \eqref{eq:3.10} we have
\begin{equation*}
\begin{split}
 & [S(t)\nu_D](x)+D^{-\frac{N}{2}}\int_0^t [S(t-s)\overline{u}(s)^q](x)\,ds\\
 & \le w(x,t)^{\frac{1}{\alpha}}+C\gamma^{\frac{pq+q}{\alpha(q+1)}-\frac{1}{\alpha}}
 [t^{-\frac{q+1}{pq-1}\alpha}]^{1-\frac{1}{\alpha}}w(x,t)^{\frac{1}{\alpha}}
 \int_0^t [s^{-\frac{q+1}{pq-1}\alpha}]^{\frac{pq+q}{\alpha(q+1)}-1}\,ds\\
 & \le w(x,t)^{\frac{1}{\alpha}}+C\gamma^{\frac{pq-1}{\alpha(q+1)}}w(x,t)^{\frac{1}{\alpha}}
\end{split}
\end{equation*}
for $(x,t)\in{\bf R}^N\times(0,1)$.
Taking a sufficiently small $\gamma>0$ if necessary,
we see that
\begin{equation}
\label{eq:3.11}
[S(t)\nu_D](x)+D^{-\frac{N}{2}}\int_0^t [S(t-s)\overline{u}(s)^q](x)\,ds\le 2w(x,t)^{\frac{1}{\alpha}}=\overline{v}(x,t)
\end{equation}
for $(x,t)\in{\bf R}^N\times(0,1)$.

Next, taking a sufficiently small $\gamma>0$ if necessary,
we show that
\begin{equation}
\label{eq:3.12}
[S(t)\mu_D](x)+D^{-\frac{N}{2}}\int_0^t [S(t-s)\overline{v}(s)^p](x)\,ds\le 2w(x,t)^{\frac{p+1}{\alpha(q+1)}}=\overline{u}(x,t)
\end{equation}
for $(x,t)\in{\bf R}^N\times(0,1)$.
We consider the case of $p\ge\alpha$.
It follows from $\alpha>1$ that
\begin{equation*}
  -\frac{q+1}{pq-1}\alpha \left( \frac{p}{\alpha}-1 \right)
  = \frac{(q+1)(\alpha-p)}{pq-1}
  > \frac{(q+1)(1-p)}{pq-1}
  = -1 + \frac{q-p}{pq-1} \ge -1.
\end{equation*}
Then, by \eqref{eq:2.1}, \eqref{eq:3.7} and \eqref{eq:3.9}
we have
\begin{equation}
\label{eq:3.13}
\begin{split}
 & [S(t)\mu_D](x)+D^{-\frac{N}{2}}\int_0^t [S(t-s)\overline{v}(s)^p](x)\,ds\\
 & \le w(x,t)^{\frac{p+1}{\alpha(q+1)}}
+C\|w(t)\|_\infty^{1-\frac{p+1}{\alpha(q+1)}}w(x,t)^{\frac{p+1}{\alpha(q+1)}}\int_0^t \|w(s)\|_\infty^{\frac{p}{\alpha}-1}\,ds\\
 & \le w(x,t)^{\frac{p+1}{\alpha(q+1)}}+C\gamma^{\frac{p}{\alpha}-\frac{p+1}{\alpha(q+1)}}
 [t^{-\frac{q+1}{pq-1}\alpha}]^{1-\frac{p+1}{\alpha(q+1)}}w(x,t)^{\frac{p+1}{\alpha(q+1)}}
 \int_0^t [s^{-\frac{q+1}{pq-1}\alpha}]^{\frac{p}{\alpha}-1}\,ds\\
 & \le w(x,t)^{\frac{p+1}{\alpha(q+1)}}+C\gamma^{\frac{p}{\alpha}-\frac{p+1}{\alpha(q+1)}}w(x,t)^{\frac{p+1}{\alpha(q+1)}}
\end{split}
\end{equation}
for $(x,t)\in{\bf R}^N\times(0,1)$.
Since $p>(p+1)/(q+1)$, taking a sufficiently small $\gamma>0$ if necessary,
we see that
$$
[S(t)\mu_D](x)+D^{-\frac{N}{2}}\int_0^t [S(t-s)\overline{u}(s)^q](x)\,ds\le 2w(x,t)^{\frac{p+1}{\alpha(q+1)}}=\overline{u}(x,t)
$$
for $(x,t)\in{\bf R}^N\times(0,1)$. Thus \eqref{eq:3.12} holds in the case of $p\ge\alpha$.

We consider the case of $p<\alpha$.
It follows from \eqref{eq:2.1} that
\begin{equation}
\label{eq:3.14}
\begin{split}
 & [S(t)\mu_D](x)+D^{-\frac{N}{2}}\int_0^t [S(t-s)\overline{v}(s)^p](x)\,ds\\
 & \le w(x,t)^{\frac{p+1}{\alpha(q+1)}}
+C\int_0^t S(t-s)w(s)^{\frac{p}{\alpha}}\,ds\\
 & \le w(x,t)^{\frac{p+1}{\alpha(q+1)}}
+C\int_0^t [S(t-s)w(s)]^{\frac{p}{\alpha}}\,ds
\le w(x,t)^{\frac{p+1}{\alpha(q+1)}}
+Ctw(x,t)^{\frac{p}{\alpha}}
\end{split}
\end{equation}
for $(x,t)\in{\bf R}^N\times(0,1)$.
Since $p>(p+1)/(q+1)$, by \eqref{eq:3.9} we have
\begin{equation}
\label{eq:3.15}
\begin{split}
tw(x,t)^{\frac{p}{\alpha}} & \le t\|w(t)\|_\infty^{\frac{p}{\alpha}-\frac{p+1}{\alpha(q+1)}}
w(x,t)^{\frac{p+1}{\alpha(q+1)}}\\
 & \le C\gamma^{\frac{p}{\alpha}-\frac{p+1}{\alpha(q+1)}}
t[t^{-\frac{q+1}{pq-1}\alpha}]^{\frac{p}{\alpha}-\frac{p+1}{\alpha(q+1)}}
w(x,t)^{\frac{p+1}{\alpha(q+1)}}\\
 & =C\gamma^{\frac{p}{\alpha}-\frac{p+1}{\alpha(q+1)}}w(x,t)^{\frac{p+1}{\alpha(q+1)}}
\end{split}
\end{equation}
for $(x,t)\in{\bf R}^N\times(0,1)$.
By \eqref{eq:3.14} and \eqref{eq:3.15},
taking a sufficiently small $\gamma>0$ if necessary,
we obtain \eqref{eq:3.12} in the case of $p<\alpha$.
Thus \eqref{eq:3.12} holds.
Combining \eqref{eq:3.8}, \eqref{eq:3.11} and \eqref{eq:3.12},
we deduce that $(\overline{u},\overline{v})$ is a supersolution to problem~(P') in ${\bf R}^N\times [0,1)$.
Thus Theorem~\ref{Theorem:3.1} follows.
$\Box$
\vspace{5pt}

\noindent
{\bf Proof of Theorem~\ref{Theorem:3.2}.}
It suffices to construct a supersolution to problem~(P') in ${\bf R}^N\times[0,1)$.
Let $L\ge e$ be such that
$\Phi_L(s):=s[\log(L+s)]^\beta$ $(s\ge 0)$ satisfies the following properties:
\begin{itemize}
  \item $\Phi_L$ is convex in $[0,\infty)$;
  \item $(0,1)\ni s\mapsto s^{\frac{p-1}{2}}[\log(L+s)]^{-p\beta}$ is monotone increasing.
\end{itemize}
Consider problem
\begin{equation}
\label{eq:3.16}
\left\{
\begin{array}{ll}
\partial_t w=\Delta w+D^{-\frac{N}{2}}w^{1+\frac{2}{N}},
& \quad x\in{\bf R}^N,\,\,t>0,\vspace{3pt}\\
w(x,0)=\displaystyle{\Phi_L^{-1}\left(\frac{1}{2}\Phi_L(2\mu_D)+\frac{1}{2}\Phi_L(2\nu_D)\right)},
 & \quad x\in{\bf R}^N.
\end{array}
\right.
\end{equation}
It follows from \eqref{eq:3.3} that
$$
\|S(t)\Phi_L(w(0))\|_\infty\le C\gamma t^{-\frac{N}{2}}\left|\log\frac{t}{2}\right|^{\beta-\frac{N}{2}}
$$
for $0<t<1$.
Taking a sufficiently small $\gamma>0$,
by the same arguments as in the proof of \cite[Theorem~5.3]{IKO}
(see also the proof of Proposition~\ref{Proposition:4.1} and  \cite[Theorem~1.5]{HI01})
we see that problem~\eqref{eq:3.16} possesses a solution~$w$ in $\mathbf{R}^N \times [0,1)$.
On the other hand,
since  $\Phi$ is convex, it follows that
$$
\Phi_L(\mu_D(x)+\nu_D(x))=\Phi_L\left(\frac{2\mu_D(x)+2\nu_D(x)}{2}\right)
\le\frac{1}{2}(\Phi_L(2\mu_D(x))+\Phi_L(2\nu_D(x)))
$$
for $x\in{\bf R}^N$.
This implies that $(w,w)$ is a supersolution to problem~(P') in $\mathbf{R}^N \times [0,1)$.
Thus Theorem~\ref{Theorem:3.2} follows.
$\Box$
\vspace{5pt}

\noindent
{\bf Proof of Theorem~\ref{Theorem:3.3}.}
Set
\begin{equation}
\notag
w(x,t):=[S(t)\mu_D](x)+[S(t)\nu_D](x).
\end{equation}
It follows from \eqref{eq:3.4} that
$\|w(t)\|_\infty\le C \gamma t^{-\frac{N}{2}}$ for $0<t<1$.
Since $q<1+2/N$, we have
\begin{equation}
\label{eq:3.17}
\begin{split}
 & [S(t)\nu_D](x) + D^{-\frac{N}{2}}\int_0^t [S(t-s)(2w(s))^q](x)\,ds\\
 & \le w(x,t)+C\int_0^t\|w(s)\|_\infty^{q-1}[S(t-s)w(s)](x)\,ds\\
 & \le w(x,t)+C\gamma^{q-1}w(x,t)\int_0^t s^{-\frac{N}{2}(q-1)}\,ds
\le w(x,t)+C\gamma^{q-1}w(x,t).
\end{split}
\end{equation}
for $(x,t)\in{\bf R}^N\times(0,1)$.
Similarly, since $p<1+2/N$, we have
\begin{equation}
\label{eq:3.18}
[S(t)\mu_D](x) + D^{-\frac{N}{2}}\int_0^t [S(t-s)(2w(s))^p](x)\,ds
\le w(x,t)+C\gamma^{p-1}w(x,t).
\end{equation}
for $(x,t)\in{\bf R}^N\times(0,1)$.
By \eqref{eq:3.17} and \eqref{eq:3.18},
taking a sufficiently small $\gamma>0$,
we see that $(2w,2w)$ is a supersolution to problem~(P') in ${\bf R}^N\times[0,1)$.
Thus Theorem~\ref{Theorem:3.3} follows.
$\Box$

\section{Case (B)}
In this section we obtain sufficient conditions on the existence of solutions
to problem~(P) in case~(B).
We prove the following theorem.
\begin{theorem}
\label{Theorem:4.1}
Let $N\ge 1$ and $0<p\le q$ with $pq>1$ be in case {\rm (B)}.
Let $\alpha>0$ and $0<\beta<1/(pq-1)$. Set
$$
\Psi(\tau):=\tau[\log(e+\tau)]^\alpha,\qquad
\Phi(\tau):=\tau [\log(e+\tau)]^\beta\quad\mbox{for}\quad\tau\ge 0.
$$
Let
$$
 \frac{q+1}{p+1}<r_*<q.
$$
Then there exists $\gamma>0$ such that,
if $\mu$ and $\nu$ are nonnegative measurable functions in ${\bf R}^N$ and satisfy
\begin{equation}
  \label{eq:4.1}
  \begin{split}
    |||S(t)\Psi(\mu)|||_{r_*}
     & \le \gamma t^{-\frac{N}{2}\left(\frac{p+1}{q+1}-\frac{1}{r_*}\right)}
    \left|\log\frac{t}{2}\right|^{-\frac{p}{pq-1}+\alpha},
   \quad 0<t<1,
    \\
    \|S(t)\Phi(\nu)\|_\infty
    & \le \gamma t^{-\frac{N}{2}}\left|\log\frac{t}{2}\right|^{-\frac{1}{pq-1}+\beta},
    \quad 0<t<1,
  \end{split}
\end{equation}
then problem~{\rm (P)} possesses a solution in ${\bf R}^N\times[0,1)$.
\end{theorem}
Similarly to Section~3,
for the proof of Theorem~\ref{Theorem:4.1},
it suffices to consider the case where $T=1$ and $D'=1$.
Let $(p,q)$ be in the case of (B).
Then
\begin{equation}
  \label{eq:4.2}
  -\frac{N}{2}(pq-1) + q = -1.
\end{equation}
Let $K_1$ and $K_2$ be positive constants such that
\begin{equation}
\label{eq:4.3}
K_1a^p+K_1b^p\ge D^{-\frac{N}{2}}(a+b)^p,
\qquad
K_2a^q+K_2b^q\ge D^{-\frac{N}{2}}(a+b)^q,
\end{equation}
for $a$, $b\ge 0$.
We obtain sufficient conditions on the existence of solutions to the Cauchy problem
\begin{equation}
\tag{Q}
\left\{
\begin{array}{ll}
\partial_t u=\Delta u+K_1v^p & \quad\mbox{in}\quad{\bf R}^N\times(0,1),\vspace{3pt}\\
\partial_t v=\Delta v+K_2u^q & \quad\mbox{in}\quad{\bf R}^N\times(0,1),\vspace{3pt}\\
u,v\ge 0 & \quad\mbox{in}\quad{\bf R}^N\times(0,1),\vspace{3pt}\\
(u(\cdot,0),v(\cdot,0))=(\mu_D,\nu_D) & \quad\mbox{in}\quad{\bf R}^N,
\end{array}
\right.
\end{equation}
under the assumption either $\mu=0$ or $\nu=0$,
and prove Theorem~\ref{Theorem:4.1}.

On the other hand, for any $L\ge e$ and $\lambda>0$, we set
$$
\Lambda(s):=s[\log(e+s)]^\lambda,\quad
\Lambda_L(s):=s[\log(L+s)]^\lambda,\quad s\ge 0.
$$
Then
\begin{equation}
\label{eq:4.4}
\begin{split}
 & C^{-1}\Lambda_L(s)\le\Lambda(s)\le C\Lambda_L(s),\\
 & 0\le\Lambda_L'(s)\le C[\log(L+s)]^\lambda,\quad
 0\le\Lambda_L'(\Lambda_L^{-1}(s))\le C[\log(L+s)]^\lambda,\\
 & C^{-1} s [\log(L+s)]^{-\lambda}
 \le \Lambda_L^{-1}(s)\le C s[\log(L+s)]^{-\lambda},
\end{split}
\end{equation}
for $s\ge 0$.
Furthermore, for any $a>0$ and $b>0$,
taking a sufficiently large $L$ if necessary, we have
\begin{itemize}
  \item[(a)] $\Lambda_L$ is convex in $[0,\infty)$;
  \item[(b)] the function $(0,1)\ni s\mapsto s^a[\log(L+s)]^{-b}$ is monotone increasing.
\end{itemize}
We prove the following proposition on the existence of solutions to problem~(Q) with $\mu_D=0$.
\begin{proposition}
\label{Proposition:4.1}
Assume the same conditions as in Theorem~{\rm\ref{Theorem:4.1}} and $\mu=0$ in ${\bf R}^N$.
Then there exists $\gamma>0$ such that,
if $\nu$ is a nonnegative measurable function in ${\bf R}^N$ satisfying
\begin{equation}
\label{eq:4.5}
\|S(t)\Phi(\nu)\|_\infty\le \gamma t^{-\frac{N}{2}}\left|\log\frac{t}{2}\right|^{-\frac{1}{pq-1}+\beta},
\quad 0<t<1,
\end{equation}
then problem~{\rm (Q)} possesses a solution in ${\bf R}^N\times[0,1)$.
\end{proposition}
{\bf Proof.}
Let $L\ge e$ and set
$$
v_*(x,t) := \Phi_L^{-1}\left[ S(t)\Phi_L (\nu_D) \right],
$$
where $\Phi_L(s) := \Lambda_L(s)$ with $\lambda = \beta$.
Let $0<\gamma<1$ and assume \eqref{eq:4.5}.
It follows from \eqref{eq:4.4} and \eqref{eq:4.5} that
\begin{equation}
\label{eq:4.6}
\|S(t)\Phi_L(\nu_D)\|_\infty\le C\gamma t^{-\frac{N}{2}}\left|\log\frac{t}{2}\right|^{-\frac{1}{pq-1}+\beta}
\equiv\gamma\xi(t),
\quad 0<t<1.
\end{equation}
\underline{Case $p>1$} :
Consider the case of $p>1$.
Taking a sufficiently large $L$ if necessary,
we can assume that properties~(a) and (b) hold with $a=(p-1)/2$ and $b=\beta p$.
Set
$$
a(t):=t^{-\frac{N}{2}(p-1)+1}\left|\log\frac{t}{2}\right|^{-\frac{p-1}{pq-1}-\beta},
\quad
U(t):=a(t)S(t)\Phi_L(\nu_D),
\quad
V(t):=2v_*(t).
$$
We show that $(U,V)$ is a supersolution to problem~(Q) in $\mathbf{R}^N \times [0,1)$.
It follows from \eqref{eq:4.2} and \eqref{eq:4.5} that
\begin{equation}
\label{eq:4.7}
\begin{split}
  & a(s)^q\|S(s)\Phi_L(\nu_D)\|_\infty^{q-1}\\
 & \le s^{-\frac{N}{2}(pq-q)+q}\left|\log \frac{s}{2}\right|^{-\frac{pq-q}{pq-1}-\beta q}\cdot
 (C\gamma)^{q-1}s^{-\frac{N}{2}(q-1)}\left|\log \frac{s}{2}\right|^{-\frac{q-1}{pq-1}+\beta(q-1)}\\
 & \le C\gamma^{q-1}s^{-1}\left|\log \frac{s}{2}\right|^{-1-\beta}
\end{split}
\end{equation}
for $0<s<1$.
Furthermore, by \eqref{eq:4.4}
we have
\begin{equation}
\label{eq:4.8}
\begin{split}
0 & \le\frac{\Phi_L(v_*(x,t))}{v_*(x,t)}=\frac{[S(t)\Phi_L(\nu_D)](x)}{\Phi_L^{-1}([S(t)\Phi_L(\nu_D)](x))}
\le C[\log(L+[S(t)\Phi_L(\nu_D)](x))]^\beta\\
 & \le C[\log(L+\gamma\xi(t))]^\beta\le C[\log(L+\xi(t))]^\beta
\le C\left|\log\frac{t}{2}\right|^\beta
\end{split}
\end{equation}
for $(x,t)\in{\bf R}^N\times(0,1)$.
By \eqref{eq:4.7} and \eqref{eq:4.8} we see that
\begin{equation}
\label{eq:4.9}
\begin{split}
 & \int_0^t [S(t-s)U(s)^q](x)\,ds
 =\int_0^t [S(t-s)a(s)^q(S(s)\Phi_L(\nu_D))^q](x)\,ds\\
 & \qquad
 \le [S(t)\Phi_L(\nu_D)](x)\int_0^t a(s)^q\|S(s)\Phi_L(\nu_D)\|_\infty^{q-1}\,ds\\
 & \qquad
 \le \left\|\frac{\Phi_L(v_*(t))}{v_*(t)}\right\|_\infty v_*(x,t)\int_0^t a(s)^q\|S(s)\Phi_L(\nu_D)\|_\infty^{q-1}\,ds\\
 & \qquad
 \le C\left|\log \frac{t}{2}\right|^\beta v_*(x,t)\int_0^t a(s)^q\|S(s)\Phi_L(\nu_D)\|_\infty^{q-1}\,ds
 \le C\gamma^{q-1}v_*(x,t)
\end{split}
\end{equation}
for $(x,t)\in{\bf R}^N\times(0,1)$.
Taking a sufficiently small $\gamma>0$ if necessary,
by \eqref{eq:2.1} and \eqref{eq:4.9} we obtain
\begin{equation}
\label{eq:4.10}
\begin{split}
 & [S(t)\nu_D](x)+K_2\int_0^t [S(t-s)U(s)^q](x)\,ds\\
 & \le v_*(x,t)+CK_2\gamma^{q-1}v_*(x,t)
\le 2v_*(x,t)=V(x,t)
\end{split}
\end{equation}
for $(x,t)\in{\bf R}^N\times(0,1)$.

On the other hand, by property~(b) (with $a=(p-1)/2$ and $b=\beta p$), \eqref{eq:4.4} and \eqref{eq:4.6} we have
\begin{equation}
\label{eq:4.11}
\begin{split}
0 & \le\frac{v_*(x,t)^p}{[S(t)\Phi_L(\nu_D)](x)}
\le C[S(t)\Phi_L(\nu_D)](x)^{p-1}[\log(L+[S(t)\Phi_L(\nu_D)](x))]^{-p\beta}\\
 & =C[S(t)\Phi_L(\nu_D)](x)^{\frac{p-1}{2}}
 [S(t)\Phi_L(\nu_D)](x)^{\frac{p-1}{2}}[\log(L+[S(t)\Phi_L(\nu_D)](x))]^{-p\beta}\\
 & \le C[\gamma\xi(t)]^{\frac{p-1}{2}}[\gamma\xi(t)]^{\frac{p-1}{2}}[\log(L+\gamma\xi(t))]^{-p\beta}\\
 & \le C[\gamma\xi(t)]^{\frac{p-1}{2}}\xi(t)^{\frac{p-1}{2}}[\log(L+\xi(t))]^{-p\beta}\\
 & \le C\gamma^{\frac{p-1}{2}}t^{-\frac{N}{2}(p-1)}\left|\log \frac{t}{2}\right|^{-\frac{p-1}{pq-1}-\beta}
\end{split}
\end{equation}
for $(x,t)\in{\bf R}^N\times(0,1)$.
Recalling that $p<1+2/N$ and taking a sufficiently small $\gamma$ if necessary,
by Lemma~\ref{Lemma:2.3} and \eqref{eq:4.11} we have
\begin{equation}
  \label{eq:4.12}
  \begin{aligned}
   & K_1\int_0^t [S(t-s)V(s)^p](x)\,ds
   = 2^pK_1\int_0^t [S(t-s)v_*(s)^p](x)\,ds\\
   & \qquad\quad
   \le 2^pK_1\int_0^t\left\|\frac{v_*(s)^p}{S(s)\Phi_L(\nu_D)}\right\|_\infty [S(t-s)S(s)\Phi_L(\nu_D)](x)\,ds\\
   & \qquad\quad
   \le C\gamma^{\frac{p-1}{2}}[S(t)\Phi(\nu_D)](x)
   \int_0^t s^{-\frac{N}{2}(p-1)}\left|\log \frac{s}{2}\right|^{-\frac{p-1}{pq-1}-\beta}\,ds\\
   & \qquad\quad
   \le C\gamma^{\frac{p-1}{2}}t^{-\frac{N}{2}(p-1)+1}\left|\log \frac{t}{2}\right|^{-\frac{p-1}{pq-1}-\beta}
   [S(t)\Phi_L(\nu_D)](x)\\
   & \qquad\quad
   =C\gamma^{\frac{p-1}{2}}a(t)[S(t)\Phi_L(\nu_D)](x)\le a(t)[S(t)\Phi_L(\nu_D)](x)=U(x,t)
 \end{aligned}
\end{equation}
for $(x,t)\in{\bf R}^N\times(0,1)$.
Therefore, by \eqref{eq:4.10} and \eqref{eq:4.12}
we see that $(U,V)$ is a supersolution to problem~(Q) in $\mathbf{R}^N \times [0,1)$
in the case of $\mu=0$.
Thus Proposition~\ref{Proposition:4.1} follows in the case of $p\ge 1$.
\vspace{5pt}
\newline
\underline{Case $0<p\le 1$} :
Consider the case of $0<p<1$.
Let $\delta\in(0,1)$ be such that
\begin{equation}
\label{eq:4.13}
\frac{N}{2}p(1-\delta)<1\quad\mbox{and}\quad\delta pq>1.
\end{equation}
Taking a sufficiently large $L$ if necessary,
we can assume that properties~(a) and (b) hold with $a=(1-\delta)/2$ and $b=\beta$.
Set
\begin{equation*}
\begin{split}
 & \tilde{a}(t):=t^{-\frac{N}{2}p(1-\delta)+1}\left|\log\frac{t}{2}\right|^{-p\frac{1-\delta}{pq-1}-\beta\delta p},\\
 & \tilde{U}(x,t):=\tilde{a}(t)[S(t)\Phi_L(\nu)](x)^{\delta p},\quad
V(x,t):=2v_*(x,t),
\end{split}
\end{equation*}
for $(x,t)\in{\bf R}^N\times(0,1)$.
We show that $(\tilde{U},V)$ is a supersolution to problem~(Q) in $\mathbf{R}^N \times [0,1)$.
It follows from \eqref{eq:4.4}, \eqref{eq:4.6} and property~(b) (with $a=(1-\delta)/2$ and $b=\beta$) that
\begin{equation}
\label{eq:4.14}
\begin{split}
\frac{v_*(x,t)}{[S(t)\Phi(\nu_D)]^\delta}
 & \le C[S(t)\Phi_L(\nu_D)](x)^{1-\delta}[\log(L+[S(t)\Phi_L(\nu_D)](x))]^{-\beta}\\
 & =C[S(t)\Phi_L(\nu_D)](x)^{\frac{1-\delta}{2}}[S(t)\Phi_L(\nu_D)](x)^{\frac{1-\delta}{2}}
 [\log(L+[S(t)\Phi_L(\nu_D)](x))]^{-\beta}\\
 & \le C(\gamma\xi(t))^{\frac{1-\delta}{2}}(\gamma\xi(t))^{\frac{1-\delta}{2}}[\log(L+\gamma\xi(t))]^{-\beta}\\
 & \le C\gamma^{\frac{1-\delta}{2}}\xi(t)^{1-\delta}[\log(L+\xi(t))]^{-\beta}\\
 & \le C\gamma^{\frac{1-\delta}{2}}
 t^{-\frac{N}{2}(1-\delta)}\left|\log\frac{t}{2}\right|^{-\frac{1-\delta}{pq-1}-\beta\delta}
\end{split}
\end{equation}
for $(x,t)\in{\bf R}^N\times(0,1)$.
Furthermore, by \eqref{eq:4.4} and \eqref{eq:4.6} we have
\begin{equation}
\label{eq:4.15}
\begin{split}
\frac{[S(t)\Phi_L(\nu_D)](x)}{v_*(x,t)} & \le C[\log(L+[S(t)\Phi_L(\nu_D)](x))]^\beta
\le C[\log(L+C\gamma\xi(t))]^\beta\\
 & \le C[\log(L+C\xi(t))]^\beta
\le C\left|\log\frac{t}{2}\right|^\beta
\end{split}
\end{equation}
for $(x,t)\in{\bf R}^N\times(0,1)$.
Taking a sufficiently small $\gamma>0$ if necessary,
by Lemma~\ref{Lemma:2.3}, \eqref{eq:2.1}, \eqref{eq:4.13} and \eqref{eq:4.15} we see that
\begin{equation}
\label{eq:4.16}
\begin{split}
 & K_1\int_0^t [S(t-s)V(s)^p](x)\,ds=2^pK_1\int_0^t [S(t-s)v_*(s)^p](x)\,ds\\
 & \qquad
 \le 2^pK_1\int_0^t \left\|\frac{v_*(s)}{[S(s)\Phi_L(\nu_D)]^\delta}\right\|_\infty^p
 [S(t-s)[S(s)\Phi_L(\nu_D)]^{\delta p}](x)\,ds\\
 & \qquad
 \le 2^pK_1\int_0^t \left\|\frac{v_*(s)}{[S(s)\Phi_L(\nu_D)]^\delta}\right\|_\infty^p
 [[S(t-s)S(s)\Phi_L(\nu_D)](x)]^{\delta p}\,ds\\
 & \qquad
 \le
 C \gamma^{\frac{p(1-\delta)}{2}}[S(t)\Phi_L(\nu_D)](x)^{\delta p}\int_0^t
 s^{-\frac{N}{2}p(1-\delta)}\left|\log\frac{s}{2}\right|^{-p\frac{1-\delta}{pq-1}-\beta\delta p}\,ds\\
 & \qquad
 \le C \gamma^{\frac{p(1-\delta)}{2}}t^{-\frac{N}{2}p(1-\delta)+1}
 \left|\log\frac{t}{2}\right|^{-p\frac{1-\delta}{pq-1}-\beta\delta p}[S(t)\Phi_L(\nu_D)](x)^{\delta p}\\
 & \qquad
 = C\gamma^{\frac{p(1-\delta)}{2}}\tilde{a}(t)[S(t)\Phi_L(\nu_D)](x)^{\delta p}
 \le \tilde{a}(t)[S(t)\Phi_L(\nu_D)](x)^{\delta p}=\tilde{U}(x,t)
\end{split}
\end{equation}
for $(x,t)\in{\bf R}^N\times(0,1)$.
On the other hand,
taking a sufficiently small $\gamma>0$ if necessary,
by \eqref{eq:4.2}, \eqref{eq:4.6}, \eqref{eq:4.13} and \eqref{eq:4.15} we have
\begin{align*}
  & K_2\int_0^t [S(t-s)\tilde{U}(s)^q](x)\,ds
  =K_2\int_0^t [S(t-s)\tilde{a}(s)^q[S(s)\Phi_L(\nu_D)]^{\delta pq}](x)\,ds\\
  & \qquad
  \le K_2[S(t)\Phi_L(\nu)](x)\int_0^t \tilde{a}(s)^q\left\|S(s)\Phi_L(\nu_D)\right\|_\infty^{\delta pq-1}\,ds\\
  & \qquad
  \le C\left|\log \frac{t}{2}\right|^\beta v_*(x,t)
  \int_0^t s^{-\frac{N}{2}pq(1-\delta)+q}\left|\log\frac{s}{2}\right|^{-pq\frac{1-\delta}{pq-1}-\beta\delta pq}\\
  & \qquad\qquad\qquad\qquad\qquad
  \times \gamma^{\delta pq-1}s^{-\frac{N}{2}(\delta pq-1)}
  \left|\log\frac{s}{2}\right|^{-\frac{\delta pq-1}{pq-1}+\beta(\delta pq-1)}\,ds\\
  & \qquad
  \le C\gamma^{\delta pq-1}\left|\log\frac{t}{2}\right|^\beta v_*(x,t)\int_0^t s^{-1}\left|\log\frac{s}{2}\right|^{-1-\beta}\,ds
  \le C\gamma^{\delta pq-1}v_*(x,t)\le v_*(x,t)
\end{align*}
for $(x,t)\in{\bf R}^N\times(0,1)$.
Therefore, by \eqref{eq:2.1} we see that
\begin{equation}
\label{eq:4.17}
[S(t)\nu_D](x)+K_2\int_0^t [S(t-s)\tilde{U}(s)^q](x)\,ds\le v_*(x,t)+v_*(x,t)=2v_*(x,t)=V(x,t)
\end{equation}
for $(x,t)\in{\bf R}^N\times(0,1)$.

By \eqref{eq:4.16} and \eqref{eq:4.17}
we see that $(\tilde{U},V)$ is a supersolution to problem~(Q) in $\mathbf{R}^N \times [0,1)$ in the case of $\mu=0$.
Thus Proposition~\ref{Proposition:4.1} follows in the case of $0<p\le 1$.
Therefore the proof of Proposition~\ref{Proposition:4.1} is complete.
$\Box$
\vspace{5pt}

Next we consider problem~(Q) in the case of $\nu=0$.
\begin{proposition}
\label{Proposition:4.2}
Assume the same conditions as in Theorem~{\rm \ref{Theorem:4.1}} and $\nu=0$ in ${\bf R}^N$.
Then there exists $\gamma>0$ such that,
if $\mu$ is a nonnegative measurable function in ${\bf R}^N$ satisfying
\begin{equation}
  \label{eq:4.18}
  |||S(t)\Psi(\mu)|||_{r_*}
  \le \gamma
  t^{-\frac{N}{2}\left(\frac{p+1}{q+1}-\frac{1}{r_*}\right)}
  \left|\log\frac{t}{2}\right|^{-\frac{p}{pq-1}+\alpha},
  \quad 0<t<1,
\end{equation}
for some $r_*\in \left((q+1)/(p+1), q \right)$,
then problem~{\rm (Q)} possesses a solution in ${\bf R}^N\times[0,1)$.
\end{proposition}
{\bf Proof.}
Let $\alpha>0$, $0<\beta<1/(pq-1)$, $\delta\in(0,1)$ and $L\ge e$.
By \eqref{eq:4.4} and \eqref{eq:4.18} we have
\begin{equation}
\label{eq:4.19}
|||S(t)\Psi_L(\mu)|||_{r_*}
\le C \gamma t^{-\frac{N}{2}\left(\frac{p+1}{q+1}-\frac{1}{r_*}\right)}
\left|\log\frac{t}{2}\right|^{-\frac{p}{pq-1}+\alpha},\quad 0<t<1,
\end{equation}
where $\Psi_L(s) := \Lambda(s)$ with $s=\alpha$.
Since $S(t)\Psi_L(\mu)=S(t/2)S(t/2)\Psi_L(\mu)$,
taking a sufficiently small $\gamma>0$ if necessary,
by \eqref{eq:2.4} and \eqref{eq:4.19} we have
\begin{equation}
\label{eq:4.20}
\begin{split}
|||S(t)\Psi_L(\mu)|||_r
 & \le Ct^{-\frac{N}{2}\left(\frac{1}{r_*}-\frac{1}{r}\right)}
\biggr|\biggr|\biggr|S\left(\frac{t}{2}\right)\Psi_L(\mu)\biggr|\biggr|\biggr|_{r_*}\\
 & \le \delta t^{-\frac{N}{2}\left(\frac{p+1}{q+1}-\frac{1}{r}\right)}\left|\log\frac{t}{2}\right|^{-\frac{p}{pq-1}+\alpha}
\end{split}
\end{equation}
for $0<t<1$ and $r_*\le r \le \infty$.
Let $r'>0$ be such that
\begin{equation}
\label{eq:4.21}
p^{-1}<r'<\frac{q+1}{p+1}.
\end{equation}
Let $\epsilon>0$ be such that $\epsilon<q$ and $p(q-\epsilon)>1$.
Let $a>0$ be such that
\begin{equation}
  \label{eq:4.22}
  \begin{split}
    & 0<a<\min\left\{q-r_*,p-1,\epsilon,\frac{p(q-\epsilon)-1}{q-\epsilon}\right\}\qquad\qquad\,\mbox{if}\quad p>1,\\
    & 0<a<\min\left\{q-r_*,pr'-1,\epsilon,r'\cdot\frac{p(q-\epsilon)-1}{q-\epsilon}\right\}\qquad\mbox{if}\quad 0<p\le 1,
  \end{split}
\end{equation}
and $b=\beta q$.
Then
\begin{equation}
  \label{eq:4.23}
  a<\epsilon<q,\qquad (q-\epsilon)(p-a)>1.
\end{equation}
Taking a sufficiently large $L\ge e$ if necessary,
we can assume that $\Psi_L$ and $\Phi_L$ have property~(a) and property~(b) holds.

Let $(U,V)$ be a solution to the Cauchy problem
\begin{equation}
\label{eq:4.24}
\left\{
\begin{array}{ll}
U_t=\Delta U+K_1\Psi_L'(\Psi_L^{-1}(U))[\Phi_L^{-1}(V)]^p,\qquad & x\in{\bf R}^N,\,\,t>0,\vspace{5pt}\\
V_t=\Delta V+K_2\Phi_L'(\Phi_L^{-1}(V))[\Psi_L^{-1}(U)]^q, & x\in{\bf R}^N,\,\,t>0,\vspace{5pt}\\
U(x,0)=\Psi_L(\mu),\,\,
V(x,0)=0, & x\in{\bf R}^N.
\end{array}
\right.
\end{equation}
Then it follows from property~(a) that
$(\tilde{u},\tilde{v}):=(\Psi_L^{-1}(U),\Phi_L^{-1}(V))$ is a supersolution to problem~\eqref{eq:P} with $\nu=0$.
Therefore it suffices to prove the existence of a solution~$(U,V)$ to Cauchy problem~\eqref{eq:4.24} in $\mathbf{R}^N \times [0,1)$.
Set $(U_0,V_0):=(S(t)\Psi_L(\mu),0)$.
Define $\{ (U_n,V_n) \}_{n=1}^\infty$ inductively by
\begin{equation}
  \label{eq:4.25}
  \begin{aligned}
    U_{n+1}(x,t)
    & :=
    S(t)\Psi(\mu)+K_1\int_0^t S(t-s)\Psi_L'(\Psi_L^{-1}(U_n(s)))[\Phi_L^{-1}(V_n(s))]^p\,ds,\\
    V_{n+1}(x,t)
    & :=K_2\int_0^t S(t-s)\Phi_L'(\Phi_L^{-1}(V_n(s)))[\Psi_L^{-1}(U_n(s))]^q\,ds,
  \end{aligned}
\end{equation}
for $x\in{\bf R}^N$ and $t>0$, where $n=0,1,2,\dots$.
Then
\begin{equation}
\label{eq:4.26}
\begin{split}
 & 0\le U_0(x,t)\le U_1(x,t)\le\cdots\le U_n(x,t)\le\cdots,\\
 & 0\le V_0(x,t)\le V_1(x,t)\le\cdots\le V_n(x,t)\le\cdots,
\end{split}
\end{equation}
for $x\in{\bf R}^N$ and $t>0$.

Let $\delta\in(0,1)$ be sufficiently small.
We show that
\begin{align}
\label{eq:4.27}
 & |||U_n(t)|||_r\le 2\delta
 t^{-\frac{N}{2}\left(\frac{p+1}{q+1}-\frac{1}{r}\right)}\left|\log\frac{t}{2}\right|^{-\frac{p}{pq-1}+\alpha}
 =:\delta\xi_r(t),\\
\label{eq:4.28}
 & |||V_n(t)|||_\ell\le \delta^{q-\epsilon}
 t^{-\frac{N}{2}\left(1-\frac{1}{\ell}\right)}\left|\log \frac{t}{2}\right|^{-\frac{1}{pq-1}+\beta}
 =:\delta^{q-\epsilon}\eta_\ell(t),
 \end{align}
for $0<t<1$, $r_*\le r\le\infty$, $1\le\ell\le\infty$ and $n=0,1,2,\dots$.
Since $V_0\equiv 0$,
by \eqref{eq:4.20} we see that \eqref{eq:4.27} and \eqref{eq:4.28} hold for $n=0$.

Assume that \eqref{eq:4.27} and \eqref{eq:4.28} hold for $n=k\in\{0,1,2,\dots\}$.
Recall that
\begin{equation}
  \label{eq:4.29}
  \frac{q+1}{pq-1} = \frac{N}{2}.
\end{equation}
By property~(b) (with $a=q-r_*-a$ and $b=q\alpha$),
\eqref{eq:4.4}, \eqref{eq:4.22}, \eqref{eq:4.27} and \eqref{eq:4.28} with $n=k$
we have
\begin{equation}
\label{eq:4.30}
\begin{split}
0 & \le\frac{\Phi_L'(\Phi_L^{-1}(V_k(x,t)))[\Psi_L^{-1}(U_k(x,t))]^q}{U_k(x,t)^{r_*}}\\
 & \le C[\log(L+V_k(x,t))]^\beta
U_k(x,t)^{q-r_*-a}U_k(x,t)^a[\log(L+U_k(x,t))]^{-q\alpha}\\
& \le C[\log(L+\delta^{q-\epsilon}\eta_\infty(t))]^\beta
(\delta\xi_\infty(t))^{q-r_*-a}(\delta\xi_\infty(t))^a[\log(L+\delta\xi_\infty(t))]^{-q\alpha}\\
 & \le C\delta^{q-r_*-a}[\log(L+\eta_\infty(t))]^\beta\xi_\infty(t)^{q-r_*}[\log(L+\xi_\infty(t))]^{-q\alpha}\\
 & \le C\delta^{q-r_*-a}\left|\log\frac{t}{2}\right|^\beta t^{-\frac{N}{2}\cdot \frac{p+1}{q+1}(q-r_*)}
 \left[\left|\log\frac{t}{2}\right|^{-\frac{p}{pq-1}+\alpha}\right]^{q-r_*}
 \left|\log\frac{t}{2}\right|^{-q\alpha}\\
 & = C\delta^{q-r_*-a}t^{-\frac{N}{2}\cdot \frac{p+1}{q+1}(q-r_*)}
 \left|\log\frac{t}{2}\right|^{\beta-\alpha r_* -\frac{p}{pq-1}(q-r_*)}
\end{split}
\end{equation}
for $(x,t)\in{\bf R}^N\times(0,1)$.
Then
\begin{equation}
\label{eq:4.31}
\begin{split}
 & \biggr|\biggr|\biggr|
 \int_{t/2}^t S(t-s)\Phi_L'(\Phi_L^{-1}(V_k(\cdot,s)))[\Psi_L^{-1}(U_k(\cdot,s))]^q\,ds
 \biggr|\biggr|\biggr|_\ell\\
 & \le C\int_{t/2}^t
 \Big|\Big|\Big|
 \Phi_L'(\Phi_L^{-1}(V_k(\cdot,s)))[\Psi_L^{-1}(U_k(\cdot,s))]^q
 \Big|\Big|\Big|_{\ell} \,ds\\
 & \le C\int_{t/2}^t
 \left\|\frac{\Phi_L'(\Phi_L^{-1}(V_k(\cdot,s)))[\Psi_L^{-1}(U_k(\cdot,s))]^q}{U_k(\cdot,s)^{r_*}}\right\|_\infty
 |||U_k(\cdot, s)^{r_*}|||_{\ell}\,ds\\
 & \le C\delta^{q-r_*-a}\int_{t/2}^t
 s^{-\frac{N}{2}\cdot \frac{p+1}{q+1}(q-r_*)}
 \left|\log\frac{s}{2}\right|^{\beta-\alpha r_* -\frac{p}{pq-1}(q-r_*)}
 |||U_k(s)|||_{\ell r_*}^{r_*}\,ds
\end{split}
\end{equation}
for $0<t<1$.
It follows from \eqref{eq:4.27} that
$$
|||U_k(s)|||_{\ell r_*}
\le 2 \delta
s^{-\frac{N}{2}\left(\frac{p+1}{q+1}-\frac{1}{\ell r_*}\right)}
\left|\log\frac{s}{2}\right|^{-\frac{p}{pq-1}+\alpha}
$$
for $0<s<1$.
This together with \eqref{eq:4.29} and \eqref{eq:4.31} implies that
\begin{equation}
\label{eq:4.32}
\begin{split}
 & \biggr|\biggr|\biggr|\int_{t/2}^t S(t-s)\Phi_L'(\Phi_L^{-1}(V_k(\cdot,s)))[\Psi_L^{-1}(U_k(\cdot,s))]^q\,ds\biggr|\biggr|\biggr|_\ell\\
 & \le C\delta^{q-a}\int_{t/2}^t
 s^{-\frac{N}{2}\left(\frac{p+1}{q+1}q-\frac{1}{\ell}\right)}
 \left|\log\frac{s}{2}\right|^{\beta-\alpha r_* -\frac{p}{pq-1}(q-r_*)
 -\frac{p}{pq-1}r_*+\alpha r_*}\,ds\\
 & \le C\delta^{q-a} \int_{t/2}^t
 s^{-\frac{N}{2}\left(1-\frac{1}{\ell}\right)-1}
 \left|\log\frac{s}{2}\right|^{\beta-1-\frac{1}{pq-1}}\,ds
 \\
 & \le C\delta^{q-a} t^{-\frac{N}{2}\left(1-\frac{1}{\ell}\right)}
 \int_{t/2}^t s^{-1}\left|\log\frac{s}{2}\right|^{\beta-1-\frac{1}{pq-1}}\,ds
\end{split}
\end{equation}
for $0<t<1$.
Similarly, by \eqref{eq:2.4} we have
\begin{equation}
\label{eq:4.33}
\begin{split}
 & \biggr|\biggr|\biggr|
 \int_0^{t/2} S(t-s)\Phi_L'(\Phi_L^{-1}(V_k(\cdot,s)))[\Psi_L^{-1}(U_k(\cdot,s))]^q\,ds\biggr|\biggr|\biggr|_\ell\\
 & \le C\int_0^{t/2} (t-s)^{-\frac{N}{2}\left(1-\frac{1}{\ell}\right)}
 |||\Phi_L'(\Phi_L^{-1}(V_k(\cdot,s)))[\Psi_L^{-1}(U_k(\cdot,s))]^q|||_1\,ds\\
 & \le C\int_0^{t/2}(t-s)^{-\frac{N}{2}\left(1-\frac{1}{\ell}\right)}
 \left\|\frac{\Phi_L'(\Phi_L^{-1}(V_k(\cdot,s)))[\Psi_L^{-1}(U_k(\cdot,s))]^q}{U_k(\cdot,s)^{r_*}}\right\|_\infty
 |||U_k(s)|||_{r_*}^{r_*}\,ds\\
 & \le C\delta^{q-a}t^{-\frac{N}{2}\left(1-\frac{1}{\ell}\right)}\int_0^{t/2}
 s^{-1}\left|
 \log\frac{s}{2}
 \right|^{\beta-\alpha r_* -\frac{p}{pq-1}(q-r_*)-\frac{p}{pq-1}r_*+\alpha r_*}\,ds\\
 & \le C\delta^{q-a}t^{-\frac{N}{2}\left(1-\frac{1}{\ell}\right)}\int_0^{t/2}s^{-1}
 \left|\log\frac{s}{2}\right|^{\beta-1-\frac{1}{pq-1}}\,ds
 \end{split}
\end{equation}
for $0<t<1$.
Since $\beta<1/(pq-1)$,
applying \eqref{eq:4.32} and \eqref{eq:4.33} to \eqref{eq:4.25},
we obtain
\begin{equation*}
  \begin{aligned}
    |||V_{k+1}(t)|||_\ell
    & \le \, K_2
    \biggr|\biggr|\biggr|
    \int_0^t S(t-s)\Phi_L'(\Phi_L^{-1}(V_k(\cdot,s)))[\Psi_L^{-1}(U_k(\cdot,s))]^q\,ds
    \biggr|\biggr|\biggr|_\ell\\
    & \le
    C\delta^{q-a}t^{-\frac{N}{2}\left(1-\frac{1}{\ell}\right)}\int_0^t s^{-1}
    \left|\log\frac{s}{2}\right|^{\beta-1-\frac{1}{pq-1}}\,ds\\
    & \le
    C\delta^{q-a}t^{-\frac{N}{2}\left(1-\frac{1}{\ell}\right)}
    \left|\log\frac{t}{2}\right|^{-\frac{1}{pq-1}+\beta}
  \end{aligned}
\end{equation*}
for $0<t<1$.
Therefore,
taking a sufficiently small $\delta>0$ if necessary,
by \eqref{eq:4.23} we obtain \eqref{eq:4.28} with $n=k+1$.

We prove \eqref{eq:4.27} with $n=k+1$.
Let us consider the case of $p>1$.
Then, by \eqref{eq:4.4}, \eqref{eq:4.22}, \eqref{eq:4.27}, \eqref{eq:4.28} and property~(b)
we have
\begin{equation}
\label{eq:4.34}
\begin{split}
0 & \le\Psi_L'(\Psi_L^{-1}(U_k(x,t)))\frac{[\Phi_L^{-1}(V_k(x,t))]^p}{V_k(x,t)}\\
 & \le C[\log(L+U_k(x,t))]^\alpha V_k(x,t)^{p-1-a}V_k(x,t)^a[\log(L+V_k(x,t))]^{-\beta p}\\
 & \le C[\log(L+\delta\xi_\infty(t))]^\alpha (\delta^{q-\epsilon}\eta_\infty(t))^{p-1-a}
 (\delta^{q-\epsilon}\eta_\infty(t))^a[\log(L+\delta^{q-\epsilon}\eta_\infty(t))]^{-\beta p}\\
 & \le C(\delta^{q-\epsilon})^{p-1-a}[\log(L+\xi_\infty(t))]^\alpha \eta_\infty(t)^{p-1}[\log(L+\eta_\infty(t))]^{-\beta p}\\
 & \le C(\delta^{q-\epsilon})^{p-1-a}
 \left[t^{-\frac{N}{2}}\biggr|\log\frac{t}{2}\biggr|^{-\frac{1}{pq-1}+\beta}\right]^{p-1}
 \biggr|\log\frac{t}{2}\biggr|^{\alpha-\beta p}
\end{split}
\end{equation}
for $(x,t)\in{\bf R}^N\times(0,1)$.
This together with \eqref{eq:4.28} implies that
\begin{equation}
\label{eq:4.35}
\begin{split}
 & \biggr|\biggr|\biggr|\int_{t/2}^t S(t-s)
 \Psi_L'(\Psi_L^{-1}(U_k(s)))[\Phi_L^{-1}(V_k(s))]^p\,ds\biggr|\biggr|\biggr|_r \\
 & \le C\int_{t/2}^t |||\Psi_L'(\Psi_L^{-1}(U_k(s)))[\Phi_L^{-1}(V_k(s))]^p|||_r\,ds\\
 & \le C(\delta^{q-\epsilon})^{p-1-a}\int_{t/2}^t
 \left[s^{-\frac{N}{2}}\biggr|\log\frac{s}{2}\biggr|^{-\frac{1}{pq-1}+\beta}\right]^{p-1}
 \biggr|\log\frac{s}{2}\biggr|^{\alpha-\beta p}||V_k(s)||_r\,ds\\
 & \le C(\delta^{q-\epsilon})^{p-a}t
 \left[t^{-\frac{N}{2}}\biggr|\log\frac{t}{2}\biggr|^{-\frac{1}{pq-1}+\beta}\right]^{p-1}
 \biggr|\log\frac{t}{2}\biggr|^{\alpha-\beta p}t^{-\frac{N}{2}\left(1-\frac{1}{r}\right)}\left|\log\frac{t}{2}\right|^{-\frac{1}{pq-1}+\beta}\\
 & \le C(\delta^{q-\epsilon})^{p-a}t^{-\frac{N}{2}\left(p-\frac{1}{r}\right)+1}\left|\log\frac{t}{2}\right|^{-\frac{p}{pq-1}+\alpha}\\
 & =C(\delta^{q-\epsilon})^{p-a} t^{-\frac{N}{2}
 \left(\frac{p+1}{q+1}-\frac{1}{r}\right)}\left|\log\frac{t}{2}\right|^{-\frac{p}{pq-1}+\alpha}
\end{split}
\end{equation}
for $0<t<1$.
Here we used the relation
\begin{equation}
  \label{eq:4.36}
  -\frac{N}{2}p + 1
  = -\frac{N}{2}p + \frac{N}{2} \cdot \frac{pq-1}{q+1}
  = -\frac{N}{2} \cdot \frac{p+1}{q+1},
\end{equation}
which follows from \eqref{eq:4.29}.
On the other hand,
since $p > 1$,
by \eqref{eq:2.4}, \eqref{eq:4.28} and \eqref{eq:4.34}
we have
\begin{equation}
\label{eq:4.37}
\begin{split}
 & \biggr|\biggr|\biggr|\int_0^{t/2} S(t-s)\Psi_L'(\Psi_L^{-1}(U_k(s)))[\Phi_L^{-1}(V_k(s))]^p \, ds\biggr|\biggr|\biggr|_r\\
 & \le C\int_0^{t/2}
 (t-s)^{-\frac{N}{2} \left(1-\frac{1}{r}\right)} |||\Psi_L'(\Psi_L^{-1}(U_k(s)))[\Phi_L^{-1}(V_k(s))]^p|||_1\,ds \\
 & \le C(\delta^{q-\epsilon})^{p-1-a}t^{-\frac{N}{2}\left(1-\frac{1}{r}\right)}
 \int_0^{t/2} s^{-\frac{N}{2}(p-1)}
 \left|\log\frac{s}{2}\right|^{-\frac{p-1}{pq-1}+\beta(p-1)+\alpha-\beta p} |||V_k(s)|||_{1}\,ds\\
 & \le C(\delta^{q-\epsilon})^{p-a}t^{-\frac{N}{2}\left(1-\frac{1}{r}\right)}\int_0^{t/2}s^{-\frac{N}{2}(p-1)}
\left|\log\frac{s}{2}\right|^{-\frac{p-1}{pq-1}-\beta+\alpha}\left|\log\frac{s}{2}\right|^{-\frac{1}{pq-1}+\beta}\,ds
\end{split}
\end{equation}
for $0<t<1$.
Recalling $p<1+2/N$ and \eqref{eq:4.36} and
combining Lemma~\ref{Lemma:2.3} and \eqref{eq:4.37},
we obtain
\begin{equation}
\label{eq:4.38}
\begin{split}
 & \biggr|\biggr|\biggr|\int_0^{t/2} S(t-s)\Psi_L'(\Psi_L^{-1}(U_k(s)))[\Phi_L^{-1}(V_k(s))]^p \, ds\biggr|\biggr|\biggr|_r \\
 & \le C(\delta^{q-\epsilon})^{p-a}t^{-\frac{N}{2}\left(1-\frac{1}{r}\right)}t^{-\frac{N}{2}(p-1)+1}
 \left|\log\frac{t}{2}\right|^{-\frac{p}{pq-1}+\alpha}\\
 & =C(\delta^{q-\epsilon})^{p-a}t^{-\frac{N}{2}\left(\frac{p+1}{q+1}-\frac{1}{r}\right)}
 \left|\log\frac{t}{2}\right|^{-\frac{p}{pq-1}+\alpha}
\end{split}
\end{equation}
for $0<t<1$.
Taking a sufficiently small $\delta$ if necessary,
by \eqref{eq:4.20}, \eqref{eq:4.23}, \eqref{eq:4.35} and \eqref{eq:4.38}
we obtain
\begin{align*}
   & |||U_{k+1}(t)|||_r\\
  & \le |||S(t)\Psi_L(\mu)|||_r  + K_1 \biggr|\biggr| \biggr|
  \int_0^t S(t-s)\Psi_L'(\Psi_L^{-1}(U_k(\cdot,s)))[\Phi_L^{-1}(V_k(\cdot,s))]^p\,ds \biggr|\biggr|\biggr|_r\\
  & \le \left[\delta+C(\delta^{q-\epsilon})^{p-a}\right]
  t^{-\frac{N}{2}\left(\frac{p+1}{q+1}-\frac{1}{r}\right)}
  \left|\log\frac{t}{2}\right|^{-\frac{p}{pq-1}+\alpha}\\
  & \le 2\delta t^{-\frac{N}{2}\left(\frac{p+1}{q+1}-\frac{1}{r}\right)}
  \left|\log\frac{t}{2}\right|^{-\frac{p}{pq-1}+\alpha}
\end{align*}
for $0<t<1$.
This implies that \eqref{eq:4.27} holds for $n=k+1$ in the case of $p>1$.

We consider the case of $0<p\le 1$.
Similarly to \eqref{eq:4.34},
by \eqref{eq:4.4}, \eqref{eq:4.22}, \eqref{eq:4.27}, \eqref{eq:4.28} with $n=k$
and property~(b) we have
\begin{equation}
\label{eq:4.39}
\begin{split}
0 & \le\frac{\left[\Psi_L'(\Psi_L^{-1}(U_k(x,t)))[\Phi_L^{-1}(V_k(x,t))]^p\right]^{r'}}{V_k(x,t)}\\
 & \le C[\log(L+U_k(x,t))]^{\alpha r'}V_k(x,t)^{pr'-1-a}V_k(x,t)^a[\log(L+V_k(x,s))]^{-\beta pr'}\\
 & \le C[\log(L+\delta\xi_\infty(t))]^{\alpha r'}(\delta^{q-\epsilon}\eta_\infty(t))^{pr'-1-a}\\
 & \qquad\times(\delta^{q-\epsilon}\eta_\infty(t))^a
 [\log(L+\delta^{q-\epsilon}\eta_\infty(t))]^{-\beta pr'}\\
 & \le C(\delta^{q-\epsilon})^{pr'-1-a}[\log(L+\xi_\infty(t))]^{\alpha r'} (\eta_\infty(t))^{pr'-1}
 [\log(L+\eta_\infty(t))]^{-\beta pr'}\\
 & \le C(\delta^{q-\epsilon})^{pr'-1-a}\left|\log\frac{t}{2}\right|^{\alpha r'}
 \left[t^{-\frac{N}{2}}\left|\log\frac{t}{2}\right|^{-\frac{1}{pq-1}+\beta}\right]^{pr'-1}
 \left|\log\frac{t}{2}\right|^{-\beta pr'}\\
 & \le C(\delta^{q-\epsilon})^{pr'-1-a}t^{-\frac{N}{2}(pr'-1)}
 \left|\log\frac{t}{2}\right|^{-\frac{pr'-1}{pq-1}+\alpha r'-\beta}
\end{split}
\end{equation}
for $(x,t)\in{\bf R}^N\times(0,1)$,
where $r'$ is as in \eqref{eq:4.21}.
By \eqref{eq:2.4}, \eqref{eq:4.28}, \eqref{eq:4.36} and \eqref{eq:4.39}
we have
\begin{equation}
\label{eq:4.40}
\begin{split}
 &  \biggr|\biggr|\biggr|\int_{t/2}^t S(t-s)
 \Psi_L'(\Psi_L^{-1}(U_k(s)))[\Phi_L^{-1}(V_k(s))]^p\,ds \biggr|\biggr|\biggr|_r \\
 & \le C\int_{t/2}^t
 \left\|\frac{\Psi_L'(\Psi_L^{-1}(U_k(s)))[\Phi_L^{-1}(V_k(s))]^p}{V_k(s)^{1/r'}}\right\|_\infty
 |||V_k(s)^{1/r'}|||_r\,ds\\
 & \le C(\delta^{q-\epsilon})^{\frac{pr'-1-a}{r'}}t^{-\frac{N}{2}\left(p-\frac{1}{r'}\right)}
 \left|\log\frac{t}{2}\right|^{-\frac{pr'-1}{r'(pq-1)}+\alpha-\frac{\beta}{r'}}
 \int_{t/2}^t |||V_k(s)|||_{r/r'}^{1/r'}\,ds\\
 & \le C(\delta^{q-\epsilon})^{\frac{pr'-1-a}{r'}}t^{-\frac{N}{2}\left(p-\frac{1}{r'}\right)}
 \left|\log\frac{t}{2}\right|^{-\frac{pr'-1}{r'(pq-1)}+\alpha-\frac{\beta}{r'}}\\
 & \qquad\qquad\times (\delta^{q-\epsilon})^{\frac{1}{r'}}t^{-\frac{N}{2}\left(\frac{1}{r'}-\frac{1}{r}\right)+1}
 \left|\log\frac{t}{2}\right|^{-\frac{1}{r'(pq-1)}+\frac{\beta}{r'}}\\
 & \le C(\delta^{q-\epsilon})^{\frac{pr'-a}{r'}}
 t^{-\frac{N}{2}\left(p-\frac{1}{r}\right)+1}\left|\log\frac{t}{2}\right|^{-\frac{p}{pq-1}+\alpha}\\
 & =C(\delta^{q-\epsilon})^{\frac{pr'-a}{r'}} t^{-\frac{N}{2}
 \left(\frac{p+1}{q+1}-\frac{1}{r}\right)}\left|\log\frac{t}{2}\right|^{-\frac{p}{pq-1}+\alpha}
\end{split}
\end{equation}
for $0<t<1$.
Furthermore,
by \eqref{eq:2.4}, Lemma~\ref{Lemma:2.3}, \eqref{eq:4.28}, \eqref{eq:4.36} and \eqref{eq:4.39}
we obtain
\begin{equation}
\label{eq:4.41}
\begin{split}
 & \biggr|\biggr|\biggr| \int_0^{t/2} S(t-s)\Psi_L'(\Psi_L^{-1}(U_k(s)))[\Phi_L^{-1}(V_k(s))]^p\,ds \biggr|\biggr|\biggr|_r\\
 & \le C\int_0^{t/2} (t-s)^{-\frac{N}{2}\left(\frac{1}{r'}-\frac{1}{r}\right)}
 |||\Psi_L'(\Psi_L^{-1}(U_k(s)))[\Phi_L^{-1}(V_k(s))]^p|||_{r'}\,ds\\
 & \le C(\delta^{q-\epsilon})^{\frac{pr'-1-a}{r'}}t^{-\frac{N}{2}\left(\frac{1}{r'}-\frac{1}{r}\right)}\int_0^{t/2}
 s^{-\frac{N}{2}\left(p-\frac{1}{r'}\right)}\left|\log\frac{s}{2}\right|^{-\frac{pr'-1}{r'(pq-1)}+\alpha-\frac{\beta}{r'}}
 |||V_k(s)|||_1^{1/r'}\,ds\\
 & \le C(\delta^{q-\epsilon})^{\frac{pr'-1-a}{r'}}t^{-\frac{N}{2}\left(\frac{1}{r'}-\frac{1}{r}\right)}\int_0^{t/2}
 s^{-\frac{N}{2}\left(p-\frac{1}{r'}\right)}\left|\log\frac{s}{2}\right|^{-\frac{pr'-1}{r'(pq-1)}+\alpha-\frac{\beta}{r'}}\\
 & \qquad\qquad\qquad\qquad\qquad\qquad\qquad\qquad\times
 (\delta^{q-\epsilon})^{\frac{1}{r'}}\left|\log\frac{s}{2}\right|^{-\frac{1}{r'(pq-1)}+\frac{\beta}{r'}}\,ds\\
 & \le C(\delta^{q-\epsilon})^{\frac{pr'-a}{r'}}t^{-\frac{N}{2}\left(\frac{1}{r'}-\frac{1}{r}\right)}\int_0^{t/2}
 s^{-\frac{N}{2}\left(p-\frac{1}{r'}\right)}\left|\log\frac{s}{2}\right|^{-\frac{p}{pq-1}+\alpha}\,ds\\
 & \le C(\delta^{q-\epsilon})^{\frac{pr'-a}{r'}}
 t^{-\frac{N}{2}\left(p-\frac{1}{r}\right)+1}\left|\log\frac{t}{2}\right|^{-\frac{p}{pq-1}+\alpha}\\
 & =C(\delta^{q-\epsilon})^{\frac{pr'-a}{r'}} t^{-\frac{N}{2}
 \left(\frac{p+1}{q+1}-\frac{1}{r}\right)}\left|\log\frac{t}{2}\right|^{-\frac{p}{pq-1}+\alpha}
\end{split}
\end{equation}
for $0<t<1$.
Here we used the relation
$$
\frac{N}{2}\left(p-\frac{1}{r'}\right)<\frac{N}{2}\left(p-\frac{p+1}{q+1}\right)
=\frac{N}{2}\frac{pq-1}{q+1}=1,
$$
which follows from \eqref{eq:4.21}.
Therefore,
applying \eqref{eq:4.20}, \eqref{eq:4.40} and \eqref{eq:4.41} to \eqref{eq:4.25},
we obtain
\begin{align*}
  & |||U_{k+1}(t)|||_r\\
  & \le |||S(t)\Psi_L(\mu)|||_r+ K_1 \biggr|\biggr|\biggr|
  \int_0^t S(t-s)\Psi_L'(\Psi_L^{-1}(U_k(\cdot,s)))[\Phi_L^{-1}(V_k(\cdot,s))]^p\,ds\biggr|\biggr|\biggr|_r\\
  & \le \left[\delta+C(\delta^{q-\epsilon})^{\frac{pr'-a}{r'}}\right]
  t^{-\frac{N}{2}\left(\frac{p+1}{q+1}-\frac{1}{r}\right)}
  \left|\log\frac{t}{2}\right|^{-\frac{p}{pq-1}+\alpha}\\
  & =\delta\left[1+C\delta^{(q-\epsilon)(p-\frac{a}{r'})-1}\right]
  t^{-\frac{N}{2}\left(\frac{p+1}{q+1}-\frac{1}{r}\right)}
  \left|\log\frac{t}{2}\right|^{-\frac{p}{pq-1}+\alpha}
\end{align*}
for $0<t<1$.
Since $(q-\epsilon)(p-a/r')>(q-\epsilon)(p-a)$,
taking a sufficiently small $\delta>0$ if necessary,
by \eqref{eq:4.23} we obtain \eqref{eq:4.27} with $n=k+1$.
Therefore we see that \eqref{eq:4.27} and \eqref{eq:4.28} hold for $n\in\{0,1,2,\dots\}$.
Then,
by \eqref{eq:4.26}, \eqref{eq:4.27} and \eqref{eq:4.28}
we see that the limit functions
$$
U(x,t):=\lim_{n\to\infty}U_n(x,t),
\qquad
V(x,t):=\lim_{n\to\infty}V_n(x,t),
$$
can be defined for $x\in{\bf R}^N$ and $t\in(0,1)$
and $(U,V)$ is a solution to problem~\eqref{eq:4.24} in ${\bf R}^N\times[0,1)$.
This implies that problem~\eqref{eq:P} possesses a supersolution in ${\bf R}^N\times[0,1)$.
Thus Proposition~\ref{Proposition:4.2} follows.
$\Box$\vspace{5pt}

Now we are ready to complete the proof of Theorem~\ref{Theorem:4.1}.
\vspace{3pt}
\newline
{\bf Proof of Theorem~\ref{Theorem:4.1}.}
Under assumptions of Theorem~\ref{Theorem:4.1},
if $\gamma$ is sufficiently small,
then problem~(Q) in the case of $\mu=0$ posses a solution $(u_1,v_1)$ in $\mathbf{R}^N \times [0,1)$.
Similarly, problem~(Q) in the case of $\nu=0$ possesses a solution $(u_2,v_2)$ in $\mathbf{R}^N \times [0,1)$.
Set $u_*=u_1+u_2$ and $v_*=v_1+v_2$.
Then, by \eqref{eq:4.3} we have
\begin{equation*}
  \left\{
  \begin{array}{ll}
    \partial_t u_* -\Delta u_*=K_1v_1^p+K_1 v_2^p\ge D^{-\frac{N}{2}}(v_1+v_2)^p=D^{-\frac{N}{2}}v_*^p
    & \mbox{in}\quad{\bf R}^N\times(0,1),\vspace{3pt}\\
    \partial_t v_* -\Delta v_*=K_2u_1^q+K_2 u_2^q\ge D^{-\frac{N}{2}}(u_1+u_2)^q=D^{-\frac{N}{2}}u_*^q
    & \mbox{in}\quad{\bf R}^N\times(0,1),\vspace{3pt}\\
    (u_*(\cdot,0),v_*(\cdot,0))=(\mu_D,\nu_D) & \mbox{in}\quad{\bf R}^N.
  \end{array}
  \right.
\end{equation*}
This means that $(u_*,v_*)$ is a supersolution to problem~(P').
Thus problem~(P) possesses a solution in $\mathbf{R}^N \times [0,1)$.
The proof of Theorem~\ref{Theorem:4.1} is complete.
$\Box$
\section{Cases (D) and (E)}
In this section we consider cases (D) and (E),
and prove the following theorem.
\begin{theorem}
\label{Theorem:5.1}
Let $N\ge 1$ and $0<p\le q$ with $pq>1$ be in case {\rm (D)} or {\rm (E)}.
Let $f$ be a positive continuous function in $(0,1)$ such that
\begin{equation}
\label{eq:5.1}
\int_0^1 \tau^{-1}f(\tau)\,d\tau<\infty.
\end{equation}
Let
\begin{equation}
\label{eq:5.2}
  \frac{Nq}{N+2} < r_* < q.
\end{equation}
Then there exists $\gamma>0$ such that,
if $\mu$ and $\nu$ are Radon measures in ${\bf R}^N$ and satisfy
\begin{equation}
  \label{eq:5.3}
  |||S(t)\mu|||_{r_*}
  \le\gamma t^{-\frac{N}{2}\left( \frac{N+2}{Nq} - \frac{1}{r_*} \right)} f(\sqrt{t})^{\frac{1}{q}},
  \quad
  \sup_{x\in{\bf R}^N}\nu(B(x,1))\le\gamma^q,
\end{equation}
for $0<t<1$,
then problem~{\rm (P)} possesses a solution in $\mathbf{R}^N \times [0,1)$.
\end{theorem}

We modify the argument in the proof of Theorem~\ref{Theorem:4.1}
to prove Theorem~\ref{Theorem:5.1}.
As in Section~3, it suffices to consider the case of $T=1$ and $D'=1$.
\vspace{5pt}
\newline
{\bf Proof of Theorem~\ref{Theorem:5.1}.}
Let $\delta\in(0,1)$ be a sufficiently small constant.
Similarly to the proof of Theorem~\ref{Theorem:4.1},
by \eqref{eq:2.4} and \eqref{eq:5.3},
taking a sufficiently small $\gamma>0$ if necessary,
we have
\begin{equation}
\label{eq:5.4}
\begin{split}
 & |||S(t)\mu|||_r
\le C\left(\frac{t}{2}\right)^{-\frac{N}{2}\left(\frac{1}{r_*}-\frac{1}{r}\right)}
\biggr|\biggr|\biggr|S\left(\frac{t}{2}\right)\mu\biggr|\biggr|\biggr|_{r_*}
\le \delta t^{\frac{N}{2r}-\frac{N+2}{2q}}\tilde{f}(\sqrt{t})^{\frac{1}{q}},\\
 & |||S(t)\nu|||_\ell\le \delta^qt^{-\frac{N}{2}\left(1-\frac{1}{\ell}\right)},
\end{split}
\end{equation}
for $0<t<1$, where $\tilde{f}(t):=f(t/2)$ and
\begin{equation}
\label{eq:5.5}
\max\left\{r_*,\frac{1}{p}\right\}\le r\le\infty,
\qquad 1\le\ell\le\infty.
\end{equation}
Set $(u_0,v_0):=(S(t)\mu,S(t)\nu)$.
Define $\{(u_n,v_n)\}_{n=1}^\infty$ inductively by
\begin{equation}
\label{eq:5.6}
\begin{split}
 & u_{n+1}(x,t):=S(t)\mu+D^{-\frac{N}{2}}\int_0^t S(t-s)v_n(s)^p\,ds,\\
 & v_{n+1}(x,t):=S(t)\nu+D^{-\frac{N}{2}}\int_0^t S(t-s)u_n(s)^q\,ds,
\end{split}
\end{equation}
for $x\in{\bf R}^N$ and $t>0$, where $n=0,1,2,\dots$.
Then
\begin{equation}
\label{eq:5.7}
\begin{split}
 & 0\le u_0(x,t)\le u_1(x,t)\le\cdots\le u_n(x,t)\le\cdots,\\
 & 0\le v_0(x,t)\le v_1(x,t)\le\cdots\le v_n(x,t)\le\cdots,
\end{split}
\end{equation}
for $x\in{\bf R}^N$ and $t > 0$.

Let $\epsilon>0$ satisfy
\begin{equation}
  \label{eq:5.8}
  \epsilon<q, \qquad p(q-\epsilon)>1.
\end{equation}
Taking a sufficiently small $\delta\in(0,1)$ if necessary,
we show that
\begin{align}
\label{eq:5.9}
 & |||u_n(t)|||_r\le \delta t^{\frac{N}{2r}-\frac{N+2}{2q}}\tilde{f}(\sqrt{t})^{\frac{1}{q}}
 + \delta t^{-\frac{N}{2}\left(p-\frac{1}{r}\right)+1},\\
\label{eq:5.10}
 & |||v_n(t)|||_\ell\le \delta^{q-\epsilon}t^{-\frac{N}{2}\left(1-\frac{1}{\ell}\right)},
\end{align}
for $0<t<1$, $n=0,1,2,\dots$, where $r$ and $\ell$ satisfy \eqref{eq:5.5}.
By \eqref{eq:5.4} we see that \eqref{eq:5.9} and \eqref{eq:5.10} hold for $n=0$.

Assume that \eqref{eq:5.9} and \eqref{eq:5.10} hold for $n=k\in\{0,1,2,\dots\}$.
It follows from \eqref{eq:2.4} that
\begin{equation}
  \label{eq:5.11}
  \biggr|\biggr|\biggr|\int^t_{t/2}S(t-s)u_k(s)^q\,ds\biggr|\biggr|\biggr|_\ell
  \le C\int_{t/2}^t|||u_k(s)^q|||_\ell\,ds
  =C\int_{t/2}^t |||u_k(s)|||_{q\ell}^q\,ds
\end{equation}
for $0<t<1$.
In cases (D) and (E), we have
\begin{equation}
  \label{eq:5.12}
  \frac{N}{2}-\frac{N}{2}pq+q+1
  = q+1-\frac{N}{2} (pq-1)
  = (pq-1)\left( \frac{q+1}{pq-1} - \frac{N}{2} \right) > 0.
\end{equation}
Since
\begin{equation}
\label{eq:5.13}
\int_0^1 s^{-1}f(\sqrt{s})\,ds=2\int_0^1\tau^{-1}f(\tau)\,d\tau<\infty,
\end{equation}
by \eqref{eq:5.1}, \eqref{eq:5.9}, \eqref{eq:5.11} and \eqref{eq:5.12} we have
\begin{equation}
  \label{eq:5.14}
  \begin{split}
    & \biggr|\biggr|\biggr|\int^t_{t/2}S(t-s)u_k(s)^q\,ds\biggr|\biggr|\biggr|_\ell\\
    & \le C\delta^q\int_{t/2}^t \left\{
    \left[s^{\frac{N}{2}\frac{1}{q\ell}-\frac{N+2}{2q}}\tilde{f}(\sqrt{s})^{\frac{1}{q}}\right]^q
    + s^{-\frac{N}{2}\left(pq-\frac{1}{\ell}\right)+q}
    \right\} ds
    \\
    &
    \le C\delta^{q}\int_{t/2}^t
    s^{\frac{N}{2\ell}-\frac{N+2}{2}}\tilde{f}(\sqrt{s})\,ds + C \delta^q t^{\frac{N}{2\ell}-\frac{N}{2}pq+q+1}
    \\
    & \le C\delta^{q}t^{-\frac{N}{2}\left(1-\frac{1}{\ell}\right)}
    \left[ \int_{t/2}^t s^{-1}\tilde{f}(\sqrt{s})\,ds + t^{\frac{N}{2} - \frac{N}{2}pq+q+1} \right]
    \le C\delta^qt^{-\frac{N}{2}\left(1-\frac{1}{\ell}\right)}
  \end{split}
\end{equation}
for $0<t<1$.
Similarly to \eqref{eq:5.14},
by \eqref{eq:2.4} we have
\begin{equation}
  \label{eq:5.15}
  \begin{split}
    & \biggr|\biggr|\biggr|\int_0^{t/2}S(t-s)u_k(s)^q\,ds\biggr|\biggr|\biggr|_\ell\\
    & \le C\int_0^{t/2} (t-s)^{-\frac{N}{2}\left(1-\frac{1}{\ell}\right)} |||u_k(s)^q|||_1\,ds
    \le Ct^{-\frac{N}{2}\left(1-\frac{1}{\ell}\right)}\int_0^{t/2} |||u_k(s)|||_q^q\,ds\\
    &
    \le C\delta^qt^{-\frac{N}{2}\left(1-\frac{1}{\ell}\right)}\int_0^{t/2}
    \left\{ \left[s^{\frac{N}{2q}-\frac{N+2}{2q}}\tilde{f}(\sqrt{s})^{\frac{1}{q}}\right]^q
    + s^{-\frac{N}{2}\left(pq-1\right)+q} \right\} ds\\
    &
    \le C\delta^{q}t^{-\frac{N}{2}\left(1-\frac{1}{\ell}\right)}
    \left[ \int_0^{t/2} s^{-1}\tilde{f}(\sqrt{s})\,ds + t^{(pq-1)\left( \frac{q+1}{pq-1} - \frac{N}{2} \right)}\right]
    \le C\delta^qt^{-\frac{N}{2}\left(1-\frac{1}{\ell}\right)}
  \end{split}
\end{equation}
for $0<t<1$.
By \eqref{eq:5.4}, \eqref{eq:5.8}, \eqref{eq:5.14} and \eqref{eq:5.15},
taking a sufficiently small $\delta>0$ if necessary,
we have
\begin{equation}
\label{eq:5.16}
  \begin{split}
    |||v_{k+1}(t)|||_\ell & \le |||S(t)\nu|||_\ell
    +D^{-\frac{N}{2}} \biggr|\biggr|\biggr|\int_0^t S(t-s)u_k(s)^q\,ds\biggr|\biggr|\biggr|_\ell\\
     & \le \delta^qt^{-\frac{N}{2}\left(1-\frac{1}{\ell}\right)}
    +C\delta^q t^{-\frac{N}{2}\left(1-\frac{1}{\ell}\right)}
    \le \delta^{q-\epsilon}t^{-\frac{N}{2}\left(1-\frac{1}{\ell}\right)}
  \end{split}
\end{equation}
for $0<t<1$.
Then \eqref{eq:5.10} holds for $n=k+1$.

On the other hand,
by \eqref{eq:2.4} and \eqref{eq:5.10} we have
\begin{equation}
\label{eq:5.17}
\begin{split}
\biggr|\biggr|\biggr|\int_{t/2}^t S(t-s)v_k(s)^p\,ds\biggr|\biggr|\biggr|_r
 & \le C\int_{t/2}^t |||v_k(s)^p|||_r\,ds
=C\int_{t/2}^t |||v_k(s)|||_{pr}^p\,ds\\
 & \le C(\delta^{q-\epsilon})^p\int_{t/2}^t s^{-\frac{N}{2}\left(p-\frac{1}{r}\right)}\,ds
\le C\delta^{p(q-\epsilon)}t^{-\frac{N}{2}\left(p-\frac{1}{r}\right)+1}
\end{split}
\end{equation}
for $0<t<1$.
Since $p<1+2/N$,
if $p\ge 1$, then by \eqref{eq:2.4} we see that
\begin{equation}
\label{eq:5.18}
\begin{split}
 & \biggr|\biggr|\biggr|\int_0^{t/2}S(t-s)v_k(s)^p\,ds\biggr|\biggr|\biggr|_r\\
 & \le C\int_0^{t/2} (t-s)^{-\frac{N}{2}\left(1-\frac{1}{r}\right)} |||v_k(s)^p|||_1\,ds
 \le Ct^{-\frac{N}{2}\left(1-\frac{1}{r}\right)}
 \int_0^{t/2} |||v_k(s)|||_p^p\,ds\\
 & \le C\delta^{p(q-\epsilon)}t^{-\frac{N}{2}\left(1-\frac{1}{r}\right)}\int_0^{t/2}
 s^{-\frac{N}{2}(p-1)}\,ds
 \le C\delta^{p(q-\epsilon)}t^{-\frac{N}{2}\left(1-\frac{1}{r}\right)-\frac{N(p-1)}{2}+1}\\
 & \le C\delta^{p(q-\epsilon)}t^{-\frac{N}{2}\left(p-\frac{1}{r}\right)+1}
\end{split}
\end{equation}
for $0<t<1$.
If $0<p<1$, then by \eqref{eq:2.4} we observe that
\begin{equation}
\label{eq:5.19}
\begin{split}
 & \biggr|\biggr|\biggr|\int_0^{t/2}S(t-s)v_k(s)^p\,ds\biggr|\biggr|\biggr|_r\\
 & \le C\int_0^{t/2}(t-s)^{-\frac{N}{2}\left(p-\frac{1}{r}\right)} |||v_k(s)^p|||_{p^{-1}}\,ds
 \le Ct^{-\frac{N}{2}\left(p-\frac{1}{r}\right)}
 \int_0^{t/2} |||v_k(s)|||_1^p\,ds\\
 & \le C\delta^{p(q-\epsilon)}t^{-\frac{N}{2}\left(p-\frac{1}{r}\right)+1}
\end{split}
\end{equation}
for $0<t<1$.
By \eqref{eq:5.4}, \eqref{eq:5.17}, \eqref{eq:5.18} and \eqref{eq:5.19} we have
\begin{equation}
\notag
\begin{split}
|||u_{k+1}(t)|||_r & \le |||S(t)\mu|||_r
+ D^{-\frac{N}{2}} \biggr|\biggr|\biggr|\int_0^t S(t-s)v_k(s)^p\,ds\biggr|\biggr|\biggr|_r\\
 & \le  \delta t^{\frac{N}{2r}-\frac{N+2}{2q}}\tilde{f}(\sqrt{t})^{\frac{1}{q}}
+C\delta^{p(q-\epsilon)}t^{-\frac{N}{2}\left(p-\frac{1}{r}\right)+1}
\end{split}
\end{equation}
for $0<t<1$.
Therefore,
taking a sufficiently small $\delta>0$ if necessary,
by \eqref{eq:5.8}, we obtain \eqref{eq:5.9} with $n=k+1$.
Thus \eqref{eq:5.9} and \eqref{eq:5.10} hold for $n\in\{0,1,2,\dots\}$.
Then,
by \eqref{eq:5.7}, \eqref{eq:5.9} and \eqref{eq:5.10}
we see that the limit functions
$$
u(x,t):=\lim_{n\to\infty}u_n(x,t),
\qquad
v(x,t):=\lim_{n\to\infty}v_n(x,t),
$$
can be defined for for $x\in{\bf R}^N$ and $t\in(0,1)$
and $(u,v)$ is a solution to problem~(P') in ${\bf R}^N\times[0,1)$.
Thus Theorem~\ref{Theorem:5.1} follows.
$\Box$
\section{Proof of Theorem~\ref{Theorem:1.2}}
As application of Theorems~\ref{Theorem:3.1}, \ref{Theorem:3.2},
\ref{Theorem:3.3}, \ref{Theorem:4.1} and \ref{Theorem:5.1},
we prove Theorem~\ref{Theorem:1.2}.
\vspace{5pt}

\noindent
\textbf{Proof of Theorem~\ref{Theorem:1.2}.}
All of the statements on the nonexistence of solutions has already been proved in~\cite{FI}
as a corollary of Theorem~\ref{Theorem:1.1}. See \cite[Corollary~1.5]{FI}.
It suffices to prove the statements on the existence of solutions.

For any $T>0$, set
\begin{equation}
\label{eq:6.1}
\mu_T(K)=T^{\frac{p+1}{pq-1}}\mu(T^{\frac{1}{2}} K),
\quad
\nu_T(K)=T^{\frac{q+1}{pq-1}}\nu(T^{\frac{1}{2}} K),
\end{equation}
for Borel sets $K$ in ${\bf R}^N$.
Then problem~(P) possesses a solution in ${\bf R}^N\times[0,T)$
if problem~(P) possesses a solution in ${\bf R}^N\times[0,1)$ with the initial data $(\mu_T,\nu_T)$.
See Remark~\ref{Remark:1.1}~(ii).
\vspace{3pt}
\newline
\underline{Case (A)} :
Let $(p,q)$ be in case (A) and $T>0$. It follows that
$$
0\le \mu_T(x)\le c_{a,1} |x|^{-\frac{2(p+1)}{pq-1}},
\quad
0\le \nu_T(x)\le c_{a,2} |x|^{-\frac{2(q+1)}{pq-1}},\quad x\in{\bf R}^N.
$$
Let $\alpha>1$ satisfy
\begin{equation*}
  \frac{2(q+1)}{pq-1}\alpha < N.
\end{equation*}
Then
\begin{equation*}
\begin{split}
 & \sup_{x\in{\bf R}^N}\int_{B(x,\sigma)}\mu_T(y)^{\frac{\alpha(q+1)}{p+1}}\,dy
+\sup_{x\in{\bf R}^N}\int_{B(x,\sigma)}\nu_T(y)^\alpha\,dy\\
 & \le c_{a,1}^{\frac{\alpha(q+1)}{p+1}}
\int_{B(0,\sigma)}|y|^{-\frac{2(p+1)}{pq-1}\frac{\alpha(q+1)}{p+1}}\,dy
+c_{a,2}^\alpha \int_{B(0,\sigma)}|y|^{-\frac{2\alpha(q+1)}{pq-1}}\,dy\\
 & \le C\big(c_{a,1}^{\frac{\alpha(q+1)}{p+1}}+c_{a,2}^\alpha\big)\sigma^{N-\frac{2\alpha(q+1)}{pq-1}}
\end{split}
\end{equation*}
for $0<\sigma\le 1$.
This together with Lemma~\ref{Lemma:2.1} implies that
$$
\big\|S(t)\mu_T^{\frac{\alpha(q+1)}{p+1}}\big\|_\infty
+ \big\|S(t)\nu_T^\alpha\big\|_\infty
\le C \big(c_{a,1}^{\frac{\alpha(q+1)}{p+1}}+c_{a,2}^\alpha\big)t^{-\frac{q+1}{pq-1}\alpha}
$$
for $0<t\le 1$.
By Theorem~\ref{Theorem:3.1} we see that
if both of $c_{a,1}$ and $c_{a,2}$ are sufficiently small,
then there exists a solution in ${\bf R}^N\times[0,1)$ to problem~(P) with the initial data $(\mu_T,\nu_T)$.
Since $T$ is arbitrary, we obtain the desired conclusion for case~(A).
\vspace{5pt}
\newline
\underline{Case (B)} :
Let $(p,q)$ be in case (B) and $r_*$ be such that
$$
\frac{q+1}{p+1}<r_*<q.
$$
Let $0< \beta <1/(pq-1)$ and set
\begin{equation}
\label{eq:6.2}
\Psi(\tau):=\tau[\log(e+\tau)]^{\frac{p}{pq-1}},\qquad
\Phi(\tau) := \tau [\log (e+\tau)]^\beta,
\end{equation}
for $\tau\ge 0$.
Let $0<c_{b,1}\le 1$. Then
\begin{equation}
\label{eq:6.3}
\begin{split}
 & 0\le\Psi(\mu(x))\\
 & \le c_{b,1}|x|^{-\frac{2(p+1)}{pq-1}}\left|\log\frac{|x|}{2}\right|^{-\frac{p}{pq-1}}
\biggr[\log\left(e+c_{b,1}|x|^{-\frac{2(p+1)}{pq-1}}\left|\log\frac{|x|}{2}\right|^{-\frac{p}{p-1}}\right)\biggr]^{\frac{p}{pq-1}}\chi_{B(0,1)}\\
 & \le c_{b,1}|x|^{-\frac{2(p+1)}{pq-1}}\left|\log\frac{|x|}{2}\right|^{-\frac{p}{pq-1}}
\biggr[\log\left(e+|x|^{-\frac{2(p+1)}{pq-1}}\left|\log\frac{|x|}{2}\right|^{-\frac{p}{p-1}}\right)\biggr]^{\frac{p}{pq-1}}\chi_{B(0,1)}\\
 & \le Cc_{b,1}|x|^{-\frac{2(p+1)}{pq-1}}\chi_{B(0,1)}
\end{split}
\end{equation}
for $x\in{\bf R}^N$.
Since $(p,q)$ is in case (B),
we see that $p<q$. Furthermore, it follows from \eqref{eq:4.2} that
$$
\frac{2(p+1)}{pq-1}=N\frac{p+1}{q+1}.
$$
These imply that $\Psi(\mu)$ belongs to the Lorentz space $L^{\frac{q+1}{p+1},\infty}({\bf R}^N)$
(see e.g. \cite{G} for the definition of Lorentz spaces).
Then we apply the Young-O'Neil inequality (see e.g. \cite[Chapter~1]{G}) to obtain
\begin{equation}
  \label{eq:6.4}
  \begin{split}
    |||S(t)\Psi(\mu)|||_{r_*} & \le\|S(t)\Psi(\mu)\|_{L^{r_*}}
    \\
    & \le C\|G(t)\|_{L^{\overline{r}}}\|\Psi(\mu)\|_{L^{r_0,\infty}}
    \le Cc_{b,1}t^{-\frac{N}{2}\left(1-\frac{1}{\overline{r}}\right)}
    \le Cc_{b,1}t^{-\frac{N}{2}\left(\frac{p+1}{q+1}-\frac{1}{r_*}\right)}
  \end{split}
\end{equation}
for $0<t<1$, where
$$
r_0:=\frac{q+1}{p+1}>1,\qquad
1+\frac{1}{r_*}=\frac{1}{\overline{r}}+\frac{1}{r_0}.
$$
On the other hand,
similarly to \eqref{eq:6.3},
we have
$$
0\le \Phi(\nu(x))\le Cc_{b,2}|x|^{-N}\left|\log\frac{|x|}{2}\right|^{-\frac{1}{pq-1}-1+\beta}\chi_{B(0,1)}(x)
$$
for $x\in{\bf R}^N$. Then, by Lemma~\ref{Lemma:2.1} we have
\begin{equation}
  \label{eq:6.5}
  \begin{aligned}
    & \| S(t) \Phi(\nu) \|_\infty
    \le
    Ct^{-\frac{N}{2}} \sup_{x\in \mathbf{R}^N} \int_{B(x,\sqrt{t})} \Phi(\nu)\, dy \\
    & \le
    Cc_{b,2} t^{-\frac{N}{2}} \int_{B(0,\sqrt{t})} |y|^{-N} \left| \log\frac{|y|}{2} \right|^{-\frac{1}{pq-1}-1+\beta}\, dy
    \le Cc_{b,2} t^{-\frac{N}{2}}\left|\log\frac{t}{2}\right|^{-\frac{1}{pq-1}+\beta}
  \end{aligned}
\end{equation}
for $0<t<1$.
By \eqref{eq:6.4} and \eqref{eq:6.5} we apply Theorem~\ref{Theorem:4.1} with \eqref{eq:6.2} to see that
problem~(P) possesses a solution in ${\bf R}^N\times[0,1)$
if $c_{b,1}$ and $c_{b,2}$ are sufficiently small.
Thus Theorem~\ref{Theorem:1.2} follows in case~(B).
\vspace{5pt}
\newline
\underline{Case (C)} :
Let $(p,q)$ be in case (C).
For $\beta>0$, set
$\Phi(\tau) := \tau [\log (e+\tau)]^\beta$ for $\tau\ge 0$.
As in case~(B), by Lemma~\ref{Lemma:2.1} we see that
\begin{equation*}
  \|S(t) \Phi(\mu)\|_\infty+\|S(t) \Phi(\nu)\|_\infty
  \le C(c_{c,1}+c_{c,2})t^{-\frac{N}{2}} \left| \log \frac{t}{2} \right|^{-\frac{N}{2}+\beta}
\end{equation*}
for $0<t<1$.
Then we apply Theorem~\ref{Theorem:3.2} to see that
problem~(P) possesses a solution if $c_{c,1}$ and $c_{c,2}$ are sufficiently small.
Thus Theorem~\ref{Theorem:1.2} follows in case~(C).
\vspace{5pt}
\newline
\underline{Case (D)} :
Let $0<T<1$ and
$$
\frac{Nq}{N+2} < r_* < q.
$$
Since $h_1$ is a increasing function on $(0,1]$,
it follows that
\begin{equation}
\label{eq:6.6}
\begin{split}
0 & \le\mu_T(x)=T^{\frac{p+1}{pq-1}}\mu(T^{\frac{1}{2}}x)
=T^{\frac{p+1}{pq-1}-\frac{N+2}{2q}}|x|^{-\frac{N+2}{q}}h_1(T^{\frac{1}{2}}|x|)\chi_{B(0,1)}(T^{\frac{1}{2}}x)\\
 & \le T^{\frac{p+1}{pq-1}-\frac{N+2}{2q}}|x|^{-\frac{N+2}{q}}h_1(|x|)
\end{split}
\end{equation}
for $x\in B(0,1)$.
By \eqref{eq:2.3} and \eqref{eq:6.6} we have
\begin{equation}
\label{eq:6.7}
\begin{split}
 & \sup_{x\in{\bf R}^N}\|S(t)\mu_T\|_{L^{r_*}(B(x,\sqrt{t}))}
\le Ct^{-\frac{N}{2}\left(1-\frac{1}{r_*}\right)}\sup_{x\in{\bf R}^N}\|\mu_T\|_{L^1(B(x,\sqrt{t}))}\\
 & \qquad
 \le Ct^{-\frac{N}{2}\left(1-\frac{1}{r_*}\right)}\int_{B(0,\sqrt{t})}\mu_T(y)\,dy\\
 & \qquad
 \le Ct^{-\frac{N}{2}\left(1-\frac{1}{r_*}\right)}T^{\frac{p+1}{pq-1}-\frac{N}{2}}
 \int_{B(0,\sqrt{t})}|z|^{-\frac{N+2}{q}}h_1(|z|)\,dz\\
 & \qquad
 \le Ct^{-\frac{N}{2}\left(1-\frac{1}{r_*}\right)}T^{\frac{p+1}{pq-1}-\frac{N}{2}}
 t^{\frac{N}{2}-\frac{N+2}{2q}}h_1(\sqrt{t})
 =CT^{\frac{p+1}{pq-1}-\frac{N}{2}}t^{-\frac{N}{2}\left(\frac{N+2}{Nq}-\frac{1}{r_*}\right)}h_1(\sqrt{t})
\end{split}
\end{equation}
for $0<t<1$.
On the other hand, since
\begin{equation}
\label{eq:6.8}
\frac{N+2}{q}r_*>\frac{N+2}{q}\frac{Nq}{N+2} =N,
\end{equation}
by \eqref{eq:6.6} we apply Lemma~\ref{Lemma:2.2} with $a=(N+2)/q$
to obtain
\begin{equation*}
  \|S(t)\mu_T\|_{L^{r_*}(B(0,1)\setminus B(0,\sqrt{t}))}
  \le CT^{\frac{p+1}{pq-1}-\frac{N+2}{2q}}
  t^{-\frac{N}{2}\left( \frac{N+2}{Nq}-\frac{1}{r_*} \right)}
  \left[h_1(t^\frac{1}{6})+t^\frac{(N+2)r_*-Nq}{4qr_*} \right]
\end{equation*}
for $0<t<1$.
This together with \eqref{eq:6.7} implies that
\begin{equation}
  \label{eq:6.9}
  \begin{split}
     & ||| S(t)\mu_T|||_{r_*}=
    \| S(t)\mu_T\|_{L^{r_*}(B(0,1))}
    \\
    & \qquad\le
    \| S(t)\mu_T\|_{L^{r_*}(B(0,\sqrt{t}))}
    + \| S(t)\mu_T\|_{L^{r_*}(B(0,1)\setminus B(0,\sqrt{t}))}
    \\
    & \qquad\le
    CT^{\frac{p+1}{pq-1}-\frac{N+2}{2q}}
    t^{-\frac{N}{2}\left( \frac{N+2}{Nq} - \frac{1}{r_*} \right)}
    \left[ \left( h_1(\sqrt{t})+h_1(t^\frac{1}{6}) + t^\frac{(N+2)r_*-Nq}{4qr_*} \right)^q \right]^\frac{1}{q}\\
    & \qquad\le
    CT^{\frac{p+1}{pq-1}-\frac{N+2}{2q}}
    t^{-\frac{N}{2}\left( \frac{N+2}{Nq} - \frac{1}{r_*} \right)}
    \left[h_1(\sqrt{t})^q+h_1(t^\frac{1}{6})^q+t^\frac{(N+2)r_*-Nq}{4r_*}\right]^\frac{1}{q}
  \end{split}
\end{equation}
for $0<t<1$.
Since
\begin{equation}
  \notag
  \int_0^1 h_1(\tau)^q \tau^{-1} \, d\tau < \infty,
\end{equation}
for any $k>0$, it follows that
\begin{equation}
  \notag
  \int_0^1 h_1(\tau^k)^q \tau^{-1} \, d\tau
  = \frac{1}{k} \int_0^1 h_1(\tau)^q \tau^{-1} \, d\tau < \infty.
\end{equation}
This together with \eqref{eq:6.8} implies that
\begin{equation}
  \label{eq:6.10}
  \int_0^1 \left[h_1(\sqrt{\tau})^q+h_1(\tau^\frac{1}{6})^q+\tau^\frac{(N+2)r_*-Nq}{4r_*}\right] \tau^{-1} \,d\tau<\infty.
\end{equation}
On the other hand, it follows from \eqref{eq:6.1} that
\begin{equation}
 \label{eq:6.11}
\sup_{x\in{\bf R}^N}\nu_T(B(x,1))=T^{\frac{q+1}{pq-1}-\frac{N}{2}}\sup_{x\in{\bf R}^N}\nu(B(x,T^{\frac{1}{2}}))
\le T^{\frac{q+1}{pq-1}-\frac{N}{2}}\sup_{x\in{\bf R}^N}\nu(B(x,1)).
\end{equation}
Since
\begin{equation*}
  \frac{p+1}{pq-1} - \frac{N+2}{2q}
  = \frac{2(q+1)-N(pq-1)}{2q(pq-1)}
  = \frac{1}{q} \left( \frac{q+1}{pq-1} - \frac{N}{2} \right) > 0,
\end{equation*}
taking a sufficiently small $T>0$ if necessary,
by \eqref{eq:6.9}, \eqref{eq:6.10} and \eqref{eq:6.11}
we apply Theorem~\ref{Theorem:5.1} to see that
problem~(P) possesses a solution in ${\bf R}^N\times[0,1)$ with the initial data $(\mu_T,\nu_T)$.
This means that problem~(P) possesses a local-in-time solution.
Thus Theorem~\ref{Theorem:1.2} follows in case~(D).
\vspace{5pt}
\newline
\underline{Case (E)} :
Let $0<T<1$ and $r_*\in (1,q)$.
It follows that
\begin{equation}
  \label{eq:6.12}
  \int_0^1f(r) r^{-1}  \, dr < \infty,
  \quad\mbox{where}\quad
  f(r) := \left[ \int_0^r h_2(\tau) \tau^{-1} \, d\tau \right]^q.
\end{equation}
Since $h_2$ is a positive function in $(0,1)$,
similarly to \eqref{eq:6.6},
we see that
\begin{equation}
  \notag
  \begin{aligned}
    \sup_{x\in \mathbf{R}^N}\int_{B(x,\sqrt{t})} \mu_T(y,0)\, dy
    & \le  C T^{\frac{p+1}{pq-1}-\frac{N}{2}} \int_{B(0,\sqrt{t})} |y|^{-N} h_2(|y|) \, dy \\
    & \le C T^{\frac{p+1}{pq-1}-\frac{N}{2}} f(\sqrt{t})^\frac{1}{q}
  \end{aligned}
\end{equation}
for $0<t<1$.
Then, by \eqref{eq:2.3} we have
\begin{equation}
\label{eq:6.13}
  \sup_{x\in \mathbf{R}^N} \|S(t)\mu_T\|_{L^{r_*}(B(x,\sqrt{t}))}
  \le C T^{\frac{p+1}{pq-1}-\frac{N}{2}} t^{-\frac{N}{2}\left( 1 - \frac{1}{r_*} \right)}f(\sqrt{t})^\frac{1}{q}
\end{equation}
for $0<t<1$.
On the other hand,
we apply Lemma~\ref{Lemma:2.2} with $a=N$ to obtain
\begin{equation*}
  \| S(t)\mu_T\|_{L^{r_*}(B(0,1)\setminus B(0,\sqrt{t}))}
  \le CT^{\frac{p+1}{pq-1}-\frac{N}{2}} t^{-\frac{N}{2}\left( 1-\frac{1}{r_*} \right)}
  \left[ h_2(t^\frac{1}{6}) + t^{\frac{N}{4}\left( 1-\frac{1}{r_*} \right)} + f(\sqrt{t})^\frac{1}{q} \right]
\end{equation*}
for $0<t<1$.
This together with \eqref{eq:6.13} implies that
\begin{equation}
  \label{eq:6.14}
  \begin{split}
    ||| S(t)\mu_T|||_{r_*}
    & =
    \| S(t)\mu_T \|_{L^{r_*}(B(0,1))}
    \\
    & \le
    \| S(t)\mu_T \|_{L^{r_*}(B(0,\sqrt{t}))}
    + \| S(t)\mu_T\|_{L^{r_*}(B(0,1)\setminus B(0,\sqrt{t}))}
    \\
    & \le
    CT^{\frac{p+1}{pq-1}-\frac{N}{2}}
    t^{-\frac{N}{2}\left( 1 - \frac{1}{r_*} \right)}
    \left[ h_2(t^\frac{1}{6})^q + t^{\frac{Nq}{4}\left( 1-\frac{1}{r_*} \right)} + f(\sqrt{t}) \right]^\frac{1}{q}
  \end{split}
\end{equation}
for $0<t<1$.
Since $q>1$, $f(1)<\infty$ and $h_2$ is increasing in $(0,1)$, we have
\begin{equation*}
  \int_0^1 h_2(\tau^\frac{1}{6})^q \tau^{-1} \, d\tau
  = 6 \int_0^1 h_2(\tau)^q \tau^{-1} \, d\tau
  \le 6h_2(1)^{q-1} f(1)^\frac{1}{q}
  < \infty.
\end{equation*}
Then, since
\begin{equation*}
  \int_0^1 \tau^{-1} f(\sqrt{\tau}) \, d\tau =
  2 \int_0^1 \tau^{-1} f(\tau) \, d\tau,
\end{equation*}
by \eqref{eq:6.12} we obtain
\begin{equation}
  \label{eq:6.15}
  \int_0^{1}
  \left[ h_2(\tau^\frac{1}{6})^q + \tau^{\frac{Nq}{4}\left( 1-\frac{1}{r_*} \right)} + f(\sqrt{\tau}) \right]
  \tau^{-1} d\tau < \infty.
\end{equation}
On the other hand,
similarly to \eqref{eq:6.11}, we have
\begin{equation}
\label{eq:6.16}
\sup_{x\in{\bf R}^N}\nu_T(B(x,1))
\le T^{\frac{q+1}{pq-1}-\frac{N}{2}}\sup_{x\in{\bf R}^N}\nu(B(x,1)).
\end{equation}
In case~(E), since $q=1+2/N>p$, it follows that
\begin{equation*}
  \frac{q+1}{pq-1} - \frac{N}{2}
  > \frac{p+1}{pq-1} - \frac{N}{2}
  = \frac{p+1}{pq-1} - \frac{1}{q-1}
  = \frac{q-p}{(pq-1)(q+1)} > 0.
\end{equation*}
Therefore, taking a sufficiently small $T>0$,
by \eqref{eq:6.14}, \eqref{eq:6.15} and \eqref{eq:6.16} we apply Theorem~\ref{Theorem:5.1} to see that
problem~(P) possesses a solution in ${\bf R}^N\times[0,1)$ with the initial data $(\mu_T,\nu_T)$.
This means that problem~(P) possesses a local-in-time solution.
Thus Theorem~\ref{Theorem:1.2} follows in case~(E).
\vspace{5pt}
\newline
\underline{Case (F)} :
By Lemma~\ref{Lemma:2.1} and \eqref{eq:6.1} we have
\begin{equation*}
  \| S(t) \mu_T \|_\infty + \| S(t) \nu_T \|_\infty
  \le C\left( T^{\frac{p+1}{pq-1} - \frac{N}{2}} + T^{\frac{q+1}{pq-1} - \frac{N}{2}} \right) t^{-\frac{N}{2}}
\end{equation*}
for $0<t<1$.
Since it follows from $p\le q$ and $q<1+\frac{2}{N}$ that
\begin{equation*}
  \frac{q+1}{pq-1} - \frac{N}{2}
  \ge \frac{p+1}{pq-1} - \frac{N}{2}
  > \frac{p+1}{pq-1} - \frac{1}{q-1}
  = \frac{q-p}{(pq-1)(q-1)} \ge 0,
\end{equation*}
taking a sufficiently small $T>0$,
we apply Theorem~\ref{Theorem:3.3} to see that
problem~(P) possesses a solution.
Thus Theorem~\ref{Theorem:1.2} follows in case~(F).
$\Box$
\bibliographystyle{amsplain}

\end{document}